\documentclass[12pt]{amsart}
\usepackage{amsmath,amscd}
\usepackage{graphicx}
\usepackage{amsfonts}
\usepackage{amssymb}

\newtheorem{thm}{Theorem}[section]
\newtheorem{lem}{Lemma}[section]
\newtheorem{prop}{Proposition}[section]
\newtheorem{cor}{Corollary}[section]

\numberwithin{equation}{section}

\newcommand{\al}{\alpha}
\newcommand{\be}{\beta}

\newcommand{\G}{\Gamma}
\newcommand{\g}{\gamma}

\newcommand{\s}{\sigma}
\newcommand{\ti}{\tilde}

\newcommand{\pa}{\partial}

\newcommand{\tH}{{\tilde{H}}}

\newcommand{\reg}{{\rm{reg}}}

\title{The mixed Hodge structure on the fundamental group of
hyperelliptic curves and higher cycles }
\author{Elisabetta Colombo}

\address{Elisabetta Colombo, Dipartimento di Matematica,
Universita' di Milano, via Saldini 50, 20133 Milano}

\email{elisabetta.colombo@mat.unimi.it}

\thanks{The  author acknowledges support from
MURST and GNSAGA (CNR) Italy}


\begin{document}

\begin{abstract}
In this paper we give a geometrical interpretation of an extension
of mixed Hodge structures (MHS) obtained from the canonical MHS on
the group ring of the fundamental group of a hyperelliptic curve
modulo the fourth power of its augmentation ideal. We show that
the class of this extension coincides with the regulator image of
a canonical higher cycle in a hyperelliptic jacobian. This higher
cycle was introduced and studied by Collino.
\end{abstract}

\maketitle
 \addtocounter{section}{-1}
\section{ Introduction}
For any pointed variety $(X,p)$ there is a canonical mixed Hodge structure (MHS) on
the group ring of its fundamental group modulo the $(k+1)$-th powers of its
augmentation ideal $J_p$. It was J.Morgan (\cite{Mor}) who first constructed such a
MHS for smooth varieties. R.Hain \cite{HD2} reformulated and extended the theory using
Chen's De Rham homotopy theory.

In the case of a pointed curve $(C,p)$ this MHS contains various
kinds of geometric information. For instance, the first of these
MHS which goes beyond cohomology, that is the MHS on $J_p/J^3_p$,
determines in general the holomorphic type of $(C,p)$ as was
proved by M.Pulte (\cite{P}, and \cite{HD1}). In the course of his
proof and using the work of B. Harris (\cite{HB}), Pulte showed
that the MHS on $J_p/J^3_p$ defines an estension class $m_p\in
Ext_{MHS}(H^3(JC,{\bf Z})_{prim},{\bf Z})\simeq
 F^2H^3(JC,{\bf C})_{prim}^*/H_3(JC,{\bf Z})_{prim}$ (where $H^*(JC,{\bf Z})_{prim}$
denotes the primitive cohomology) such that $2m_p$ equals the
Abel-Jacobi
 image in the primitive intermediate
jacobian of the so-called {\it Ceresa cycle} $C_p-C_p^-\subset
JC$.

In this paper we seek geometrical interpretations of the MHS on $J_p/J_p^4$. In
general the MHS on it has length 2 and it is an extension of $J_p/J_p^3$. To start
this investigation we restrict ourselves to hyperelliptic curves. In fact for a
hyperelliptic curve with Weierstra{\ss} point $p$ the short exact sequence
$0\rightarrow J^2_p/J^3_p\rightarrow J_p/J_p^3 \rightarrow J_p/J_p^2\rightarrow 0$
splits as extension of MHS mod torsion, corresponding to the fact that in this case
$C_p=C^-_p$ (cf. Prop.\ref{split}). Hence there is a natural way to transform parts of
the MHS on $J_p/J_p^4$ into extensions of MHS's that are purely built up from the
cohomology. As in Pulte's theorem on $J_p/J_p^3$, this allows us to compare these
extensions with known geometrical data.  Here we give an interpretation  of one of
these extensions in terms of the regulator image $reg(Z)$, where $Z$ is a canonical
higher cycle $Z$ constructed by Collino (\cite{Col}) on the jacobian of a
hyperelliptic curve with two fixed Weierstra{\ss} points $q_1$ and $q_2$. It is not
difficult to check that the other subextensions can be built from $J_p/J_p^3$. As
explained in \cite{Col} 1.1, $reg(Z)\in J_2(JC)_{prim}:=F^1H^2(JC,{\bf
C})_{prim}^*/H_2(JC,{\bf Z})_{prim}$ can be interpreted as a degeneration of the
Abel-Jacobi image of the Ceresa cycle for the singular curve obtained from $C$ by
glueing $q_1$ and $q_2$. The extension class $Pe\in Ext_{MHS}(H^2(JC)_{prim},{\bf
Z})\simeq J_2(JC)_{prim}$ that we compare with $reg(Z)$ is in fact obtained from the
MHS on the group ring of the fundamental group of the punctured curves $C-\{q_1\}$ and
$C-\{q_2\}$, modulo the fourth power of the augmentation ideal. The main result (Th.
\ref{risultato}) is the equality: $$Pe=(2g+1)reg(Z).$$

\

In section 1 we recall briefly the definition of regulator map and
we write the regulator image of the higher cycle $Z$ in terms of
integrals and iterated integrals on the curve. In section 2 we
construct the extension $Pe$ and state Theorem \ref{risultato}. In
section 3, using properties of iterated integrals and extension
theory, we prove the theorem. The proof is a very explicit
computation of both terms of the equality on a basis of
$F^1H^2(JC,{\bf C})_{prim}$. In section 4 and 5 we extend these
constructions to families. In particular
 in section 4
we provide an alternative proof of the result of Collino that $reg(Z)$ is not zero for
general hyperelliptic $C$, by showing that the homomorphism of fundamental groups
induced by the normal function extending $reg(Z)$ is not trivial (Cor.\ref{torsione}).
Notice that the non triviality of $reg(Z)$ implies that $Z$ is an indecomposable
cycle. Finally in section 5 we study the homomorphisms of fundamental groups induced
by normal functions associated to the extensions of Sec.2.

\

{\small {\it Acknowledgments} : I wish to thank A.Collino, R.Kaenders and E.Looijenga
for very useful discussions and suggestions.}

\section{The regulator map for hyperelliptic jacobians}
The main goal of this section is to write down the regulator image of the higher cycle
$Z$ constructed by Collino on the jacobian of an hyperelliptic curve in terms of
(iterated) integrals on the curve. First we recall the definition of the regulator and
the construction of $Z$.

\subsection{The regulator.}\label{Hcy}
 Let $X$ be a smooth projective variety of
dimension $n$. Let $CH^{n}(X,1)$ be the first higher Chow group. (For the general
theory of higher Chow groups and regulator maps we refer to \cite{Bl} and \cite{Be}).
An element in $CH^{n}(X,1)$ is defined by the 'cycle' $A:=\sum_i(C_i,f_i)$, where
$C_i$ is an irreducible curve on $X$ and $f_i$ is a rational function on $C_i$, such
that $\sum_i[div(f_i)]=0$ on $X$ (see \cite{F}, 1.3 pag.10 for the definition of the
0-cycle $[div(f_i)]$). Denote by $[0,\infty]$ the positive real axis on ${\bf P}^1$
and let $\gamma_i:=\mu_{i*}f_i^{-1}([0,\infty])$, where $\mu_{i}:
\tilde{C_i}\rightarrow C_i$ is the resolution of singularities of the curve $C_i$ and
$f_i$ is viewed as a map $f_i:\tilde{C_i}\rightarrow {\bf P}^1$. The condition
$\sum_i[div(f_i)]=0$ implies that the 1-chain $\sum \g_i$ is a 1-cycle and in fact up
to torsion it is a boundary. Suppose $H_1(X,{\bf Z})$ has no torsion (otherwise the
construction has to be done over ${\bf Q}$), the  Bloch-Beilinson regulator map in
this particular case is defined in the following way:
\begin{equation}
\label{def1reg} \reg : CH^n(X,1)\rightarrow J_2(X):=
{\frac{(F^1H^2(X,{\bf C}))^*}{H_2(X,{\bf Z}(1))}}, \quad
A=\sum_i(C_i,f_i)\mapsto \reg(A)
\end{equation}
where $\reg(A)$ is the class of the current:
\begin{equation}\label{reg}
\al\mapsto\sum_i\int_{C_i-singC_i}log(f_i) \al+ 2\pi i\int_D\al
\end{equation}
with $\al$ a closed 2-form on $X$ whose
 cohomology class is in
$F^1H^2(X,{\bf C})$ and $D$ a 2-chain such that $\pa D=\sum \g_i$. An analogous
definition can be given using instead the primitive cohomology:
\begin{equation}\label{def2reg}
reg:CH^n(X,1)\longrightarrow
 J_2(X)_{prim}:=
{{F^1H^{2}(X,{\bf C})_{prim}^*}\over{H_2(X,{\bf Z}(1))_{prim}}}
\end{equation}

\subsection{The Collino cycle}\label{Z}
Our variety $X$ will be the jacobian $JC$ of a hyperelliptic curve $C$ of genus $g$.
We will follow \cite{Col} to construct a canonical higher cycle,
$$Z:=(C_1,h_1)+(C_2,h_2),$$ as follows. For Weierstra{\ss} points $q_1$ and $q_2$, let
$h$ be a 2 to 1 morphism
\begin{equation}\label{h}
h:C\longrightarrow {\bf P}^1,\quad h(q_1)=0,\; h(q_2)=\infty.
\end{equation}
Call $C_s$ (for $s=1,2$) the image of $C$ via the
Abel Jacobi map:
\begin{equation}
\mu_{q_s}:C\longrightarrow JC,\quad x\mapsto x-{q_s},\quad C_s:=\mu_{q_s}(C).
\end{equation}
If we think of $\mu_{q_s}$ as a biholomorphism onto its image we
can define $h_s:C_s\longrightarrow {\bf P}^1,\quad h_s:=h\circ
\mu_{q_s}^{-1}$. The rational functions $h_1$ and $h_2$  on $C_1$
and $C_2$ satisfy $div (h_1)=2\sigma-2O=-div (h_2)$, where $O$ is
the origin of $JC$ and $\sigma:=q_2-q_1\in JC$. The point $\sigma$
is
 2 torsion, since $q_1$ and $q_2$ are Weierstra{\ss} points. Thus $Z=(C_1,h_1)+(C_2,h_2)\in
 CH^g(JC,1)$ and
 Collino proved that $reg(Z)$ is not zero, for
general $C$.

\subsection{Basic properties of iterated integrals}

We recall the definition of iterated integrals. Suppose $X$ is a
smooth manifold. Given $\g:[0,1]\rightarrow X$ a piece-wise smooth
path in $X$, and smooth 1 forms $w_1,...w_l$ on $X$ we define an
iterated integral by: $$\int_{\g}w_1...w_l:=\int...\int_{0\leq
t_1\leq...\leq t_l\leq 1} f_1(t_1)f_2(t_2)...f_l(t_l)dt_1...dt_l$$
where $\g^*w_j=f_j(t)dt$.

For an introduction to the subject we refer to \cite{HD1}. Here we
only collect some properties of iterated integral that we need in
the sequel and that are easy to check from the definition:

\begin{lem}  If $\omega_1$ and $\omega_2$ are
smooth 1-forms, $df$ is a smooth exact 1-form and  $\al$ and
$\beta$ are piece-wise smooth paths in $X$ with $\al(1)=\be(0)$
then: $$
\begin{array}{rl}1.& \int_{\al\cdot
\be}\omega_1\omega_2=\int_{\al}\omega_1\omega_2+
\int_{\be}\omega_1\omega_2+\int_{\al}\omega_1 \int_{\be}\omega_2\\
2.&
\int_{\al}\omega_1\omega_2+\int_{\al}\omega_2\omega_1=\int_{\al}\omega_1
\int_{\al}\omega_2\\
 3.& \int_{\al}df\omega_1=
\int_{\al}f\omega_1- f(\al(0))\int_{\al}\omega_1\\
 4.& \int_{\al}\omega_1df=
f(\al(1))\int_{\al}\omega_1-\int_{\al}f\omega_1 .
 \end{array}
$$\label{ii}
\end{lem}

\subsection{} In the next lemmas we show that for $\reg(Z)$ the expression $\int_D$ in
\ref{reg} can be written as an iterated integral. First we fix
some notations. Let $[0,\infty]$ be the positive real axis in
${\bf P}^1$, $\g:=h^{-1}([0,\infty])$, with $h$ as in \ref{h}.
Since $h$ is a covering of degree 2 outside the set of
Weierstra{\ss} points we can write $\g$ as a union of two paths
$\g =\g^++\g^-$ where $\g^+$ and $\g^-$ live in different sheets
of $h$ and have in common just the Weierstra{\ss} points in $\g$.
We fix the diffeomorphism $[0,1] \rightarrow [0,\infty]\subset
{\bf P}^1$, $t\mapsto {t\over{1-t}}$ and a $C^{\infty}$
parametrization of $\g^{\pm}$ compatible with it, i.e.:
$$\g^{\pm}:[0,1]\longrightarrow C\quad t\mapsto \g^{\pm}(t)\quad
{\rm with}\quad h(\g^{\pm}(t))={t\over{1-t}}\in [0,\infty]\subset
{\bf P}^1$$

\subsection{Remark} \label{gammas}We also consider, for $s=1,2$,
$\g_s=\mu_{q_s}(\g)$, written as $\g_s=\g_s^++\g_s^-$ with
$\g_s^{\pm}=\mu_{q_s}(\g^{\pm}([0,1]))$. Hence $\g_s^{\pm}$ are
paths with
 $C^{\infty}$ parameterizations $\mu_s\circ \g^{\pm}$.

 In the next lemma we show that
$\g_1+\g_2$ is the boundary of a 2 chain $D$ in
 $JC$ which is the sum of two
parameterized disks $D^+$ and $D^-$.

\begin{lem}Let $D^{\pm}\subset JC$ be the images of
$$F^{\pm}:[0,1]\times [0,1]\longrightarrow JC,\quad (t,s)\mapsto \g^{\pm}\left
(1-{{t(1-s)}\over{1-s(1-t)}}\right)-\g^{\pm}(1-t),$$ then $D^{\pm}$ has boundary
$\partial D^{\pm}$ with parameterization $\gamma_1^{\pm}\cdot\g_2^{\mp}$ and hence
$\partial D=\partial (D^++D^-)=\gamma_1+\g_2.$
\end{lem}
\begin{proof} Restrict the map $F^{+}$ to the boundary of
 $[0,1]\times [0,1]$:

$$
\begin{array}{lr}
s=0,&\quad \{\g^+(1-t)-\g^+(1-t)\}=O\\ t=1,&\quad
\{\g^+(s)-q_1\}=\g_1^+\subset C_1\\
 s=1,&\quad
\{q_2-\g^+(1-t)\}=\{\g^-(1-t)-q_2\}=(\g_2^-)^{-1}\subset C_2\\
t=0,&\quad \{q_2-q_2\}=O
\end{array}$$
so the oriented boundary of $D^+$ is $\partial
D^+=\g_1^+\cdot\g_2^-$. (Note that $\g_1^+(0)=O=\g_2^-(1)$ and
 $\g_1^+(1)=q_2-q_1=\s=\g_2^-(0)$.)
 The same computation yields  $\partial
D^-=\g_1^-\cdot\g_2^+$

\end{proof}

\begin{lem} Let $\phi,\psi$ be closed 1-forms with $\psi$ of type 1,0
on $J(C)$. Then,  for $D^{\pm}$ as in the previous lemma,
$$\int_{D^{\pm}}\phi\wedge \psi=\int_{\g_1^{\pm}}\phi\psi-
\int_{\g_2^{\mp}}\psi\phi,$$ where $\int\phi\psi$ denotes the
iterated  integral.\label{D+-}
\end{lem}
\begin{proof} Every closed form on a disk is exact, so
$\phi_{|D^{\pm}}=d\rho^{\pm}$, then by Stokes theorem: $$\int_{D^{\pm}}\phi\wedge
\psi=\int_{D^{\pm}}d(\rho^{\pm}\psi)= \int_{\partial D^{\pm}}\rho^{\pm}
\psi=\int_{\g_1^{\pm}\cdot\g_2^{\mp}} \rho^{\pm} \psi .$$ Moreover by
Lemma.\ref{ii}.(3), choosing $\rho^{\pm}(0)=0$: $$\int_{D^{\pm}}\phi\wedge
\psi=\int_{\g_1^{\pm}\cdot\g_2^{\mp}} d\rho^{\pm}
\psi=\int_{\g_1^{\pm}\cdot\g_2^{\mp}} \phi \psi$$ and by Lemma.\ref{ii}.(1)
$$=\int_{\g_1^{\pm}} \phi \psi+\int_{\g_2^{\mp}} \phi \psi+\int_{\g_1^{\pm}} \phi
\int_{\g_2^{\mp}} \psi.$$ Now note that $\int_{\g_1^{\pm}} \phi+\int_{\g_2^{\mp}}
\phi=0$ since $\g_1^{\pm}\g_2^{\mp}=\partial D^{\pm}$ is homotopically trivial, so
$$\int_{D^{\pm}}\phi\wedge \psi=\int_{\g_1^{\pm}} \phi \psi+ \int_{\g_2^{\mp}} \phi
\psi-\int_{\g_2^{\mp}} \phi\int_{\g_2^{\mp}} \psi=\int_{\g_1^{\pm}} \phi \psi-
\int_{\g_2^{\mp}} \psi \phi$$ where the last equality follows from Lemma.\ref{ii}(2).

\end{proof}

\

 We are ready to compute $\reg(Z)$. It is enough to do it on harmonic
forms:

\begin{thm} Let $\phi$ and $\psi$ be
harmonic 1-forms on $J(C)$ with $\psi$ of type 1,0 and denote in
the same way the corresponding 1-forms on $C$ .  Then:
$$\reg(Z)(\phi\wedge \psi)=2\int_{C-\g}log(h)\phi\wedge \psi+ 2\pi
i \int_{\g} (\phi \psi- \psi \phi).$$ \label{regreg}
\end{thm}

\begin{proof} From Lemma.\ref{D+-} $\int_{D}\phi\wedge
\psi=\int_{\g_1}\phi\psi- \int_{\g_2}\psi\phi$. Now use the fact
that the harmonic forms on $JC$ are translation invariant.
\end{proof}

\subsection{Remark}\label{closed} The right hand side of the equality of
Thm.\ref{regreg} can be computed more generally  for $\phi$ and
$\psi$ closed 1-form, since it is zero if one of them is exact.
\section{Extensions}
The theory of iterated integrals for pointed Riemann surfaces
$(C,p)$ and pointed punctured ones $(C-\{q\},p)$ describes
explicitely the canonical MHS on the quotients $J_p/J_p^k$ and
$J_{q,p}/J^k_{q,p}$, where $J_{p}:={\rm ker}(\epsilon:{\bf
Z}\pi_1(C,p)\rightarrow {\bf Z})$ and $J_{q,p}:={\rm
ker}(\epsilon:{\bf Z}\pi_1(C-\{q\},p)\rightarrow {\bf Z})$ (we
refer to \cite{HD1} and in particular for punctured curves to
\cite{Ka1}). The weight filtrations on the duals are given, for
$l\leq k$, by $W_l(J_{p}/J^k_{p})^*:=(J_{p}/J^{l+1}_{p})^*$ and
$W_l(J_{q,p}/J^k_{q,p})^*:=(J_{q,p}/J^{l+1}_{q,p})^*$. The graded
factors have the following identifications:

\begin{equation}(J^{l}_{p}/J^{l+1}_{p})^*\simeq
Q_l(C):=\bigcap_{i=0}^{i=l-2} (\otimes^iH^1(C,{\bf Z}))\otimes
Q_2(C)\otimes (\otimes^{l-2-i}H^1(C,{\bf Z}))
\end{equation}
where $Q_2:={\rm ker}(\cup:H^1(C,{\bf Z})\otimes H^1(C,{\bf
Z})\rightarrow H^2(C,{\bf Z}))$, and
\begin{equation}(J^{l}_{q,p}/J^{l+1}_{q,p})^*\simeq
\otimes^lH^1(C-\{q\},{\bf Z}) \simeq \otimes^lH^1(C,{\bf Z}).
\end{equation}The Hodge filtration is given by Chen's $\pi_1$-De Rham
-Theorem (cf.\cite{HD1}). This MHS is compatible with the natural
extensions:
\begin{equation}\label{hkp}
h^k_{p}:\quad\quad 0\rightarrow (J_{p}/J_{p}^{k-1})^*\rightarrow
(J_{p}/J_{p}^k)^*\rightarrow Q_{k-1}(C)\rightarrow 0.
\end{equation}
\begin{equation}\label{hk}
h^k_{q,p}:\quad\quad 0\rightarrow
(J_{q,p}/J_{q,p}^{k-1})^*\rightarrow
(J_{q,p}/J_{q,p}^k)^*\rightarrow \otimes^{k-1} H^1(C,{\bf
Z})\rightarrow 0.
\end{equation}

\subsection{Hyperelliptic case}\label{notazione}

Let now $C$ be a hyperelliptic curve with hyperelliptic involution
$i$. It holds:

\begin{prop} Let $C$ be a hyperelliptic
curve and let $p$ and $q$ be Weierstra{\ss} points. Then the
extensions classes $h^3_p\in Ext_{MHS}(Q_2(C),H^1(C,{\bf Z}))$ and
$h^3_{p,q}\in Ext_{MHS}(\otimes^2H^1(C,{\bf Z}),H^1(C,{\bf Z}))$
are 2-torsion, i.e. $2h^3_p=0$ and $2h^3_{q,p}=0$. \label{split}
\end{prop}
\begin{proof}
The hyperelliptic involution $i$ induces an automorphism of
$(J_{p}/J_{p}^3)^*$ such  that $h^3_p\sim i(h^3_p)=-h^3_p$. Hence
$2h^3_p=0$ in $ Ext_{MHS}(Q_2,H^1)$. The same argument holds for
$h^3_{q,p}$.

\end{proof}

 Let $p$, $q_1$ and $q_2$ be Weierstra{\ss} points on
$C$. Fix a 2 to 1 map: $h:C\rightarrow {\bf P}^1$ such that $h(q_1)=0$ and
$h(q_2)=\infty$ as in the construction of the Collino cycle $Z$.  Moreover fix a set
of loops $\{\al_l\}$ ($l=1,...,2g$) on $C$ with basepoint $p$ whose homotopy classes
give a system of generators of $\pi_1(C,p)$ with the relation:
$\prod_k[\al_k,\al_{g+k}]$. We choose all $\al_l$ not passing through $q_1$ and $q_2$
so they define also a system of generators of the  groups $\pi_1(C-\{q_s\},p)$,
($s=1,2$). Let $\{A_l\}$ ($A_l :=[\al_l],\; \in H_1(C,{\bf Z})$) be the associated
symplectic basis and $\{dx_l\}$ the dual basis of $H^1(C,{\bf Z})$. We will identify
$dx_l$ with the corresponding harmonic 1-form. From now on we denote $H^1(C,{\bf Z})$
by $H^1$.

\begin{prop} The linear map
$ r_{3,2}^s:(J_{q_s,p}/J^3_{q_s,p})^{*}\rightarrow H^1 $, which is
the dual of the linear map defined by: $A_l\mapsto
(\al_l-i_*\al_l)\; mod (J^3_{q_s,p}) $,
 is a morphism of  MHS.\label{retrazione}
 \end{prop}
\begin{proof} The space $H_{1}(C,{\bf Q})$ can be identified with the
eigenspace of the eigenvalue $-1$ of the involution $i_*\in
Aut((J_{q_s,p}/J^3_{q_s,p})_{{\bf Q}})$ ).

\end{proof}

In order to define the extension class $Pe$ of the main theorem,
we need also the following natural morphisms of MHS:
\begin{enumerate}\label{maps}

\item The monomorphism given by tensoring with the polarization
$\Omega$:
\begin{equation}\label{Omega}
J_{\Omega}=\otimes\Omega:H^1(-1)\longrightarrow \otimes^3H^1 \end{equation}
where
$H^1(-1)$ denotes $H^1$ twisted by the Tate Hodge structure ${\bf Z}(-1)$.

\item The surjection: $$\Pi:\otimes^2 H^1\rightarrow {\bf Z}(-1),$$
given by the cup product $\otimes^2H^1\rightarrow H^2$ composed
with the isomorphism $H^2(C,{\bf Z})\simeq {\bf Z}(-1)$ defined by
the integration over $C$.

\item The monomorphism $\iota:\wedge^2H^1\rightarrow H^1\otimes
H^1,\quad\quad \phi\wedge \psi\mapsto \phi\otimes \psi -
\psi\otimes \phi.$
\end{enumerate}
Notice that the map $\Pi\circ \iota$ can be identified with  the integration
$\int_C:\wedge^2 H^1\rightarrow {\bf Z}$ over $C$ and it's in this form that we shall
often write it. Set $$\wedge^2H^1_{prim}:=ker\Pi\circ \iota=ker \int_C.$$

\subsection{Main result
}\label{ext} Let $$e_s\in Ext_{MHS}(\otimes^3H^1,H^1)$$ (for
$s=1,2$) be the extension class obtained by pushing forward
$h^4_{q_s,p}$ along  $r_{3,2}^s$. The pull back of $e_2-e_1$ along
$J_{\Omega}$ defines the extension class $$e_{\Omega}\in
Ext_{MHS}(H^1(-1),H^1).$$ By tensoring  by $H^1$ on the left and
then by pushing it down along $\Pi$ one obtains $$\tilde {e}\in
Ext_{MHS}(\otimes^2 H^1(-1),{\bf Z}(-1))\simeq Ext_{MHS}(\otimes^2
H^1,{\bf Z}).$$ The pull back  along the monomorphism
$\iota:\wedge^2H^1\rightarrow \otimes^2H^1$ defines
\begin{equation}\label{e}
e\in Ext_{MHS}(\wedge^2 H^1,{\bf Z}). \end{equation} Finally the
pullback along $Ker\int_C\hookrightarrow \wedge^2H^1$
 defines $$Pe\in Ext_{MHS}(Ker\int_C,{\bf Z}).$$

The isomorphism $\wedge^2H^1\simeq H^2(JC,{\bf Z})$ and the standard theory of
separated extensions of MHS (see \cite{Ca}), that we briefly recall in the next
section, tell us that we can identify: $$ Ext_{MHS}(\wedge^2 H^1,{\bf Z})\simeq
{{Hom(H^2(JC,{\bf C}),{\bf C})}\over{F^0+ Hom(H^2(JC,{\bf Z}),{\bf Z})}}\simeq
J_2(JC).
 $$
Moreover the identification $\wedge^2H^1_{prim}$ with $ H^2(JC)_{prim}$ gives the
isomorphism $$ Ext_{MHS}(\wedge^2H^1_{prim},{\bf Z})\simeq J_2(JC)_{prim}.$$

The main result is: \begin{thm} Let $C$ be a hyperelliptic curve
and let $q_1,q_2$ and $p$ be Weierstra{\ss} points. Let
$h:C\rightarrow {\bf P}^1$ a 2:1 map with $h(q_1)=0$ and
$h(q_2)=\infty$, then
\begin{equation}
e=(2g+1)\left(\reg(Z)+log(h(p))\int_C\right)\in J_2(JC),
\end{equation}
which implies
\begin{equation}Pe=(2g+1)reg(Z)\in J_2(JC)_{prim}.
\end{equation}
\label{risultato}
\end{thm}

\section{Carlson's representatives and proof of \ref{risultato}}
To a MHS $V$ whose weights are all negative, can be associated the
intermediate jacobian:
\begin{equation}\label{JV}
 J(V):={{V_{\bf C}}\over
{F^0V_{\bf C}+V_{\bf Z}}}.
\end{equation}
An extension of MHS $0\rightarrow A\rightarrow H \rightarrow
B\rightarrow 0$ is called separated if the minimal non zero weight
of $B$ is bigger then the maximal non zero weight of $A$. Thus in
particular $Hom(B,A)$ has all negative weights. Carlson's theory
of separated extensions of MHS defines the isomorphism (see
\cite{Ca}) $$Ext_{MHS}(B,A)\simeq J(Hom(B,A)),$$ which associates
to the class of $0\rightarrow A\rightarrow H\rightarrow B
\rightarrow 0 $  the class of the composed map $ r_{{\bf Z}}\circ
s_F\in Hom(B_{{\bf C}},A_{{\bf C}})$ (called the Carlson
representative), where $s_F\in Hom(B_{{\bf C}},H_{{\bf C}})$ is a
  section preserving the Hodge filtration and $r_{{\bf Z}}\in Hom(H,A)$
is a retraction of ${\bf Z}$ modules. In this section we first
describe explicitly the Carlson representatives of the extensions
classes introduced in Sect.2, then we manipulate these expressions
using just basic properties of integrals over Riemann surfaces and
iterated integrals and at the end we prove Thm.\ref{risultato}.
\subsection{The extension class $h^4_{q_s,p}\in Ext_{MHS}(\otimes^3 H^1,
(J_{q_s,p}/J^3_{q_s,p})^*$).} Via Carlson theory, $h^4_{q_s,p}$,
defined in \ref{hk}, corresponds to the class of $(r_{4,3}^s)\circ
s^4_{Fs}$ defined in the following way. The linear map
$$r_{4,3}^s:(J_{q_s,p}/J_{q_s,p}^4)^* \rightarrow
(J_{q_s,p}/J_{q_s,p}^3)^*$$ is the dual of a linear map defined by
fixing a basis in $J_{q_s,p}/J_{q_s,p}^3$ and lifting the elements
of this basis to independent elements in $J_{q_s,p}/J_{q_s,p}^4$.
The only condition required is that the chosen basis contains the
elements $(\al_l-1)\; mod\;J^3_{q_s,p}$ which are lifted to
$(\al_l-1)\; mod\;J^4_{q_s,p}$.

 The section preserving the Hodge filtration $s^{4s}_{F}: \otimes^3
H^1\rightarrow (J_{q_s,p}/J_{q_s,p}^4)^*$ is provided by Chen
theory:

\begin{equation}
s^{4s}_{F}(dx_l\otimes dx_m \otimes dx_{n})= \int dx_ldx_mdx_{n}
+dx_l\mu_{mn,q_s} +\mu_{lm,q_s}dx_{n}+\mu_{lmn,q_s}
\end{equation}
where $\mu_{lm,q_s}$ $\mu_{mn,q_s}$ and $\mu_{lmn,q_s}$  are
smooth 1,0 logarithmic forms on $C-\{q_s\}$ satisfying:

\begin{equation} \label{A1}
dx_l\wedge dx_m +d\mu_{lm,q_s}=0, \quad dx_m\wedge dx_n
+d\mu_{mn,q_s}=0,
\end{equation}
\begin{equation} \label{A2}
dx_l\wedge\mu_{mn,q_s}+\mu_{lm,q_s}\wedge dx_{n}+
d\mu_{lmn,q_s}=0.
\end{equation}

\subsection{The extension class $e_s\in Ext_{MHS}(\otimes^3H^1,H^1)$
($s=1,2$)} Since $e_s$ is obtained by pushing forward along
$r_{3,2}^s$, it can be identified as the class of the map
$$G_s:=(r_{3,2}^s\circ r_{4,3}^s)\circ s^{4s}_F\in
Hom(\otimes^3H^1,H^1).$$ Notice that the map $r_{3,2}^s\circ
r_{4,3}^s:(J_{q_s,p}/J_{q_s,p}^4)^*\rightarrow H^1$ is in fact the
dual of the linear map given by $A_l\mapsto (\al_l-i_*\al_l) mod
(J_{q_s,p}^4)$ and it is not longer a morphism of MHS as it is
$r^s_{3,2}$.

\subsection{The extension class $e_{\Omega}\in Ext_{MHS}(H^1(-1),H^1)$}

The extension  $e_{\Omega}$ is obtained by pulling back $e_2-e_1$
along $J_{\Omega}$, hence it is representated by:
$$G:=(G_2-G_1
)\circ J_{\Omega}\in Hom(H^1(-1),H^1).$$

\subsection{}In the next proposition we compute  $G$
explicitly on the basis chosen in \ref{notazione}. Recall that
$\Omega$ can be written in coordinates as
$$\Omega=\sum_k(dx_k\otimes dx_{g+k}- dx_{g+k}\otimes dx_k)\in
\otimes^2 H^1.$$
\subsection{Assumption}\label{scelte}
Notice that we can choose solutions of \ref{A1} and \ref{A2}
satisfying the properties listed below (for $s=1,2$):
\begin{enumerate}
\item $\mu_{ml,q_s}=-\mu_{lm,q_s}$;

\item for $|l-m|\neq g$, $\mu_{lm,q_s}$ is smooth on $C$ and orthogonal to all
harmonic forms, i.e. $\mu_{lm,q_s}\wedge dx_n$ is exact;

\item $\mu_{i(g+i),q_2}$
has logarithmic singularity on $q_2$ with residue 1;

\item for $|l-m|\neq g$,
$\mu_{lm,q_1}=\mu_{lm,q_2}$;

\item $\mu_{i(g+i),q_1}=\mu_{i(g+i),q_2}+
dh/2h.$ Notice that $\mu_{i(g+i),q_1}$ has a pole of order 1 on $q_1$ with residue 1;

\item $i^*\mu_{lm,q_s}=\mu_{lm,q_s}$; $i^*\mu_{lmn,q_s}=-\mu_{lmn,q_s}$.

\end{enumerate}

\begin{prop} With the choices done
 in
Assumption \ref{scelte}, a map $G:H^1(-1)\rightarrow H^1$, whose class defines
$e_{\Omega}$, is given by:
\begin{equation}\label{G(dxl)}
G(dx_l)(A_m)=\int_{\al_m}\left[(2g+1)(log(h)-log(h(p)))dx_l
-2W(dx_l)\right]
\end{equation}
where
\begin{equation}
W(dx_l):=\sum_{k=1}^g\{(\mu_{lk(g+k),q_2}- \mu_{lk(g+k),q_1})
-(\mu_{l(g+k)k,q_2}-\mu_{l(g+k)k,q_1})\}.
\end{equation}
\label{G}
\end{prop}
\begin{proof} From the definition of $J_{\Omega}$ and $s^{4s}_{F}$ and by the
equality $\mu_{(g+k)kq_s}=-\mu_{k(g+k)q_s}$ of Assump\-tion \ref{scelte}
(1):$$G(dx_l)(A_m)=\int_{(\al_m-i_*\al_m)}
\left[\sum_k\{2dx_l(\mu_{k(g+k)q_2}-\mu_{k(g+k)q_1})\right.$$ $$\left.+(\mu_{lk,q_2}-
\mu_{lk,q_1})dx_{g+k}-(\mu_{l(g+k),q_2}- \mu_{l(g+k),q_1})dx_k\} +W(dx_l)\right].$$

By Assumption \ref{scelte}(4), setting $i=l$ if $l\leq g$ and $i=l-g$ if $l>g$:
$$G(dx_l)(A_m)=\int_{(\al_m-i_*\al_m)}
\left[2dx_l\sum_k(\mu_{k(g+k)q_2}-\mu_{k(g+k)q_1})\right. $$
$$\left.-(\mu_{i(g+i),q_2}-\mu_{i(g+i),q_1})dx_{l}+W(dx_l)\right]$$ and by Assumption
\ref{scelte}(5): $$G(dx_l)(A_m)=-\int_{(\al_m-i_*\al_m)} \left[g
dx_l(dh/h)-{{1}\over{2}}(dh/h)dx_l+W(dx_l)\right].$$ From the equalities: $$
\int_{i_*\al_m}[gdx_l(dh/h)-{{1}\over{2}}(dh/h)dx_l]=
-\int_{\al_m}[gdx_l(dh/h)-{{1}\over{2}}(dh/h)dx_l]$$ coming from $i^*dx_l=-dx_l$,
$i^*dh/h=dh/h$ and, $$ \int_{i_*\al_m}W(dx_l)=-\int_{\al_m}W(dx_l),$$ coming from
 $i^*W(dx_l)=-W(dx_l)$ (see Assumption \ref{scelte}(6))
 it follows that:

\begin{equation}
G(dx_l)(A_m)=\int_{\al_m}-[2gdx_l(dh/h)-(dh/h)dx_l -2W(dx_l)].
\end{equation}
On $C-\gamma$, $dh/h$ is an exact form. Hence if $\alpha_m\cap \gamma=\phi$, then the
statement follows from Lemma\ref{ii}(3) and (4). If $\alpha_m\cap \gamma\neq\phi$,
then the computation has to be done on a path lifting $\alpha_m$ on a covering of $C$
where $dh/h$ is exact. But the difference between it and the expression \ref{G(dxl)}
is given by a multiple of $2\pi i\int_{\alpha_m}dx_l$. Hence it defines an element in
$Hom_{\mathbf{Z}}(H^1(-1),H^1)$, which is trivial in $J(Hom(H^1(-1),H^1)).$
\end{proof}

\subsection{The extension class $\tilde{e}\in Ext_{MHS}(\otimes^2 H^1,{\bf Z})$}
The extension $\tilde{e}$ is constructed by tensoring $e_2-e_1$ by
$H^1$
 on the left and pushing forward along $\Pi$ so it can be
identified with the class of the map:
\begin{equation}\label{defF}
F:=\Pi\circ(id\times G) \in (\otimes^2 H^1)^*.
\end{equation}
\begin{prop} Let $F$ be the map defined in \ref{defF}. We have
\begin{equation} F(dx_m\otimes
dx_{l})=\end{equation} $$c(m)\int_{\al_{\sigma(m)}}\left[
(2g+1)(log(h)-log(h(p)))dx_l-2W(dx_l)\right]$$ where
$\sigma(m)=g+m$ and $c(m)=1$ if $m\leq g$, $\sigma(m)=m-g$ and
$c(m)=-1$ if $m>g$. \label{F}
\end{prop}
\begin{proof} The map $\Pi:\otimes^2 H^1\rightarrow {\bf Z}$
(cf.\ref{maps}(2)) can be written, with our choice of basis, as:
$$\Pi(v\otimes w )= \sum_{k=1}^g [v(A_k)w(A_{g+k})-
v(A_{g+k})w(A_k)].$$ Hence
 $$F(dx_m\otimes dx_l)=\Pi\circ (id\otimes G )(dx_m\otimes dx_l)=$$

$$\sum_k[dx_m(A_k)G(dx_l)(A_{g+k})-dx_m(A_{g+k})G(dx_l)(A_k)]$$ $$
=c(m)G(dx_l)(A_{\sigma(m)})$$ and finally, by \ref{G(dxl)} and by Lemma \ref{ii}(3)
$$=c(m)\int_{\al_{\sigma(m)}}(2g+1)(log(h)- \log(h(p)))dx_l-2W(dx_l).$$

\end{proof}

The next proposition is a fundamental step toward the equality of
Th.3.1. It provides an identification of the key integrals over
$C$ with iterated integrals along paths. Let
$\g:=h^{-1}([0,\infty])$ as in Sect.1.
\begin{prop}
Let $\al$ be a simple smooth loop on $C$ transverse to $\g$. Let $\phi$, $\psi$ and
$\varpi$ be 1-forms such that $\phi$, $\psi$ and $(log(h)\psi+\varpi)$ are closed and
the cohomology  class of $\phi$ is the Poincar\`e dual of $[\al]$. Then:
$$\int_{\al}(log(h)\psi+\varpi)= \int_{C-\g}\phi\wedge (log(h)\psi+\varpi)+2\pi
i\int_{\g}\phi\psi.$$\label{kra}
\end{prop}

\begin{proof} Denote by $\eta$ a closed 1-form in the same cohomology
class of $\phi$ with compact support on a tubular neighborhood  of
$\al$. Hence $\phi=\eta+df$ where $df$ is an exact 1-form. Thus
$df\wedge (log(h)\psi+\varpi) =d(f(log(h)\psi+\varpi))$ is an
exact form on $C-\g$ or, equivalently, on the Riemann surface with
boundary obtained from $C$ cutting along $\g$. Hence by Stokes
theorem and by taking the difference of the determinations for
$log$ at the boundary, we have:
\begin{equation}\label{dfm}
\int_{C-\g}df\wedge(log(h)\psi+\varpi)=-2\pi i \int_{\g}f \psi.
\end{equation}
Choose $f$ such that $f(p_0)=0$, with $p_0$ the basis point  of
$\al$, so

$$\int_{\g}f \psi=\int_{\g}(f-f(p_0))\psi=\int_{\g}df \psi.$$ To
compute $\int_C\eta \wedge (log(h)\psi+\varpi)$ we recall that the
class of $\eta$ is  the Poincar\'e dual of the class $\al$. We
give a more explicit construction of such $\eta$ with support on a
tubular neighborhood $D=D^+\cup D^-$ of $\al$. Following for
example \cite{FK}II.3.3 let $G$ be a ${\bf C}^{\infty}$ function
on $C-\al$, which is the constant $1$ on a smaller strip
$D_{0}\subset D^-$ and $0$ on $C-D^-$. Then take $\eta$ equal to
$dG$ in $D-\al$ and 0 otherwise. We distinguish two cases.
 First
suppose that $\alpha$ doesn't intersect $\gamma$. Then we can take
$D\cap \gamma=\phi$, so $log(h)\psi+\varpi$ is a closed form well
defined on the support of $\eta$, and:
\begin{equation}\label{etam1}\int_{C-\g}\eta \wedge
\left(log(h)\psi+\varpi\right)=\int_{\al}
\left(log(h)\psi+\varpi\right).\end{equation} Moreover since in this case
$\phi_{|\g}=df_{|\g}$, adding \ref{dfm} and \ref{etam1} gives the result. Suppose now
that $\al$ intersects $\g$. Notice that now $log(h)$ is not well defined on $D$ and we
need to compute the integral on disjoint union of  rectangles $D"$ obtained by cutting
$D$ along $D\cap \g$. Applying Stokes' theorem:

\begin{equation}\label{etam2}\int_{C-\g}\eta \wedge
\left(log(h)\psi+\varpi\right)=\int_{D"}
\eta\wedge\left(log(h)\psi+\varpi\right)\end{equation}

$$=\int_{\al}\left(log(h)\psi+\varpi\right)-2\pi i \int_{\g\cap
D^- }dG\psi,$$ since $G$ is $0$ outside $D_{0}$, and $\eta=dG$ on
$C-\al$, we obtain
$$=\int_{\al}\left(log(h)\psi+\varpi\right)-2\pi i \int_{\g
}\eta\psi.$$ To conclude we add the equalities \ref{dfm} and
\ref{etam2} recalling that
$\phi_{|\gamma}=\eta_{|\gamma}+df_{|\gamma}$.

\end{proof}

\begin{cor} Choosing as $\al_m$ simple smooth loops transverse to $\g$,
we get: $$ \int_{C-\g}dx_m\wedge
\left(log(h)dx_l-{{2}\over{2g+1}}W(dx_l)\right)+2\pi i\int_{\g}
dx_mdx_l= $$
\begin{equation}c(m)\int_{\al_{\sigma(m)}}\left(log(h)dx_l-{{2}\over{2g+1}}W(dx_l)\right)
\end{equation} where $\sigma(m)=g+m$ and $c(m)=1$ if $m\leq g$,
$\sigma(m)=m-g$ and $c(m)=-1$ if $m>g$.\label{integrale}
\end{cor}

\begin{proof} The 1-form $log(h)dx_l-{{2}\over{2g+1}}W(dx_l)$ is
closed and the class of $dx_m$ is $c(m)$ times the Poincar\`e dual
of the class of $\al_{\sigma(m)}$.
\end{proof}

\begin{cor}\label{integrale2}
$$ F(dx_m\otimes dx_l)=(2g+1)\left[ \int_{C-\gamma}dx_m\wedge \left(
(log(h)-log(h(p)))dx_l-{{2}\over{2g+1}}W(dx_l)\right)+2\pi
i\int_{\g}dx_mdx_l\right].$$
\end{cor}

\begin{proof} This follows from Prop.\ref{F} and
Cor.\ref{integrale}.

\end{proof}

Since we have the isomorphism
\begin{equation}
J((\otimes^2H^1)^*)\simeq {{F^1(\otimes^2H^1_{\bf
C})^*}\over{(\otimes^2 H^1)^*}},
\end{equation}
in order to determine $F$ it is enough to compute it on elements
in $F^1(\otimes ^2H^1_{\bf C})$, namely linear combinations of
$dx_l\otimes dz_i$ and $dz_i\otimes dx_l$, where
$\{dz_i\}_{i=1,...g}$ is a basis of $H^{1,0}(X)$. We choose such a
basis to satisfy the condition $\int_{\al_i}dz_j=\delta_{ij}$.
Hence $$dz_i=dx_i+\sum_{j=1}^gZ_{ij}dx_{g+j}\quad\quad{\rm
with}\quad \quad Z_{ij}:= \int_{\al_{g+i}}dz_j.$$
\begin{prop}
The map $F$ evaluated on elements $dz_i\otimes dx_l$ gives: $$F(dz_i\otimes
dx_l)=(2g+1)\left[\int_C(log(h)-log(h(p)) dz_i\wedge dx_l+2\pi
i\int_{\g}dz_idx_l\right].$$\label{dxdz}
\end{prop}
\begin{proof} This follows from  Cor. \ref{integrale2}, using the linearity
of $F$ and the fact that $dz_i\wedge W(dx_l)=0$ for reasons of
type.

 \end{proof}

In order to compute $F(dx_l\otimes dz_i)$ we prove the following:

\begin{lem} With a
suitable choices of the $\mu_{lmn,q_s}$, we have: $$W(dz_i):=W(dx_i)+\sum_{j=1}^g
Z_{ij}W(dx_{g+j})=0.$$ \label{dimW(dz)}
\end{lem}

\begin{proof}
 We set
$$\mu_{lk(g+k),q_s}=R_{lk(g+k),q_s}+S_{lk(g+k),q_s},\quad \quad
\mu_{l(g+k)k,q_s}=R_{l(g+k)k,q_s}+S_{l(g+k)k,q_s},$$ where the $R_{lk(g+k),q_s}$ and
the $R_{l(g+k)k,q_s}$ satisfy
\begin{equation} \label{D}
dx_l\wedge\mu_{k(g+k),q_s}+dR_{lk(g+k),q_s}=0,\quad\quad
dx_l\wedge\mu_{(g+k)k,q_s}+dR_{l(g+k)k,q_s}=0
\end{equation}
 while the $S_{lk(g+k),q_s}$ and the $S_{l(g+k)k,q_s}$ satisfy
\begin{equation}\mu_{lk,q_s}\wedge dx_{g+k}+dS_{lk(g+k),q_s}=0, \quad\quad
\mu_{l(g+k),q_s}\wedge dx_{k}+dS_{l(g+k)k,q_s}=0.
\end{equation}
Thus $W(dx_i)$ becomes: $$W(dx_l)=WR(dx_l)+WS(dx_l)$$ with $$WR(dx_l)=
\sum_k\{(R_{lk(g+k),q_2}-R_{lk(g+k),q_1}) -(R_{l(g+k)k,q_2}-R_{l(g+k)k,q_1})\},$$
$$WS(dx_l)=\sum_k\{(S_{lk(g+k),q_2}-S_{lk(g+k),q_1})
-(S_{l(g+k)k,q_2}-S_{l(g+k)k,q_1})\}.$$ Now we claim that:
\subsection{Claim}\label{Rz} {\it After a suitable choice of $R_{...}$ and $S_{...}$,
$$1)\quad\quad WR(dz_i):=WR(dx_i)+\sum_jZ_{ij}WR(dx_{g+j})=0$$ $$2)\quad\quad
WS(dz_i):=WS(dx_i)+\sum_jZ_{ij}WS(dx_{g+j})=0.$$} The Lemma follows directly by
\ref{Rz}, which we prove below
\end{proof}

\begin{proof} (of Claim \ref{Rz}) Point 1): We choose
$R_{l(g+k)k,q_s}:=-R_{lk(g+k),q_s},$
 so
$$WR(dx_l)= 2\sum_k\{(R_{lk(g+k),q_2}-R_{lk(g+k),q_1}).$$ The
thesis follows once we fix any $R_{(g+j)k(g+k),q_s}$ satisfying
condition \ref{D} and define: $$R_{ik(g+k),q_s}:=-\sum_jZ_{ij}
R_{(g+j)k(g+k),q_s}.$$

Point 2): Since by Assumption.\ref{scelte} for $|m-k|\neq g$,
$\mu_{mk,q_2}=\mu_{mk,q_1}$, we can choose $S_{mk(g+k),q_2}=S_{mk(g+k),q_1}$ and for
the same argument for $|m-(g+k)|\neq g$, $S_{m(g+k)k,q_2}=S_{m(g+k)k,q_1}.$ Thus
$$WS(dx_i)=S_{i(g+i)i,q_1}-S_{i(g+i)i,q_2}\quad
WS(dx_{g+j})=S_{(g+j)j(g+j),q_2}-S_{(g+j)j(g+j),q_1}.$$
 For all $j\leq g$, fix
$S_{(g+j)j(g+j),q_s}$ ($s=1,2$) and $S_{j(g+j)j,q_1}$. To get the
result it enough to set: $$S_{i(g+i)i,q_2}:=S_{i(g+i)i,q_1}+
\sum_jZ_{ij}(S_{(g+j)j(g+j),q_2}-S_{(g+j)j(g+j),q_1}).$$

\end{proof}

\begin{prop}
The map $F$ evaluated on elements $dx_l\otimes dz_i $ gives:
$$1)\quad \quad F(dx_l\otimes
dz_i)=(2g+1)c(l)\int_{\al_{\sigma(l)}}((log(h)-log(h(p)))dz_l$$
where $\sigma(l)=g+l$ and $c(l)=1$ if $l\leq g$, $\sigma(l)=l-g$
and $c(l)=-1$ if $l>g$, $$2)\quad \quad F(dx_l\otimes
dz_i)=(2g+1)\left[\int_C(log(h)-log(h(p))) dx_l\wedge dz_i-2\pi
i\int_{\g}dx_ldz_i\right].$$ \label{dzdx} \end{prop}
\begin{proof}
First notice that
 by Prop.\ref{F} and using the linearity of
$F$ we have that $$F(dx_l\otimes
dz_i)=(2g+1)c(l)\int_{\al_{\sigma(l)}}\left[(log(h)-log(h(p)))dz_i-2W(dz_i)\right].$$
Then Point 1) follows directly from Prop.\ref{F} (using the linearity of $F$) and
Lemma \ref{dimW(dz)}. Point 2) follows from:
$$c(l)\int_{\al_{\sigma(l)}}log(h)dz_i=\int_Cdx_l\wedge log(h)dz_i+2\pi
i\int_{\g}dx_ldz_i$$  which is the equality of Prop.\ref{kra} with $\phi=c(l)dx_l$,
$\psi=dz_i$ and $\varpi=0$.

\end{proof}

\subsection{The extension classes ${e\in Ext_{MHS}(\wedge^2H^1,{\bf Z})}$ and
$Pe\in Ext_{MHS}(\wedge^2H^1_{prim},{\bf Z})$}

Carlson representatives of $e$ and $Pe$ are simply $F\circ\iota \in (\wedge^2H^1)^*$
and  $F\circ\iota_{|\wedge^2H^1_{prim}} \in (\wedge^2H^1_{prim})^*$. The final step
towards the proof of Th.\ref{risultato} is the following proposition in which we
compute $F\circ\iota$ on the elements $dx_l\wedge dz_i$ of the basis of
$F^1\wedge^2H^1$.
\begin{prop} The map $F\circ\iota$ whose class defines ${e\in Ext_{MHS}(\wedge^2H^1,{\bf
Z})}$  is given by:
\begin{equation}F\circ\iota (dx_l\wedge dz_i)=
F(dx_l\otimes dz_i)-F(dz_i\otimes dx_l)
\end{equation}
$$=(2g+1)\left[2\int_Clog(h) dx_l\wedge dz_i+2\pi
i\int_{\g}(dx_ldz_i-dz_idx_l)+2log(h(p))\int_Cdx_l\wedge dz_i\right].$$ \label{Fiota}
\end{prop}

\begin{proof} This follows from Prop.\ref{dxdz} and Prop.\ref{dzdx}.

\end{proof}

\

Now the main result follows immediately:
\begin{proof} (of Thm.\ref{risultato}) We compare the explicit expressions
in Prop.\ref{Fiota} and Thm.\ref{regreg}.

\end{proof}

\subsection{Remark} Notice that using Thm.\ref{regreg} and Prop.\ref{kra}
$\reg(Z)$ can be computed in a simple way in terms of iterated
integrals along paths on $C$, namely we have: $$\reg(Z)(dx_m\wedge
dz_i)=2c(m)\left[\int_{\al_{\sigma(m)}}(dh/h)dz_i+log(h(p))
\int_{\al_{\sigma(m)}}dz_i\right],$$ where $\sigma(m)=g+m$ and
$c(m)=1$ if $m\leq g$, $\sigma(m)=m-g$ and $c(m)=-1$ if $m>g$.

\section{The normal function defined by the regulator}

In this section we extend the construction of $reg(Z)$ to families. In this setting we
construct a normal function on a fine moduli space of hyperelliptic curves with
Weierstra{\ss} points. We show that the homomorphism between fundamental groups
induced by such a normal function is not trivial. This provides an alternative proof
(Cor.\ref{torsione}) of the result of Collino that $reg(Z)$ is not zero for general
hyperelliptic curves. The method of proof of Collino was to show that the associated
infinitesimal invariant of a normal function extending $reg(Z)$ was not zero. For all
the theory related to moduli spaces of curves and mapping class groups we refer to
\cite{HL}.

\subsection{Mapping class group}\label{S}
 Fix a compact
orientable surface $S$ of genus $g$ together with $n$ distinct
points $x_1,...,x_n$. The mapping class group $\G^n_{g}$ is the
group of isotopy classes of orientation preserving diffeomorphisms
of $S$ that fix each of the chosen points (for $n=0$ we will drop
the apex). A classically known system of generators of
 $\G^n_{g}$  is given by the Dehn twists  $D_a$ of
simple closed curves $a\subset S$. The mapping class group
$\G^n_{g}$ has a natural representation
\begin{equation}\label{ro1}
\rho:\G^n_{g}\longrightarrow Sp(H_1(S,\mathbf{Z}))\simeq Sp_{g}(\mathbf{Z})
\end{equation}
given by the action of $\G^n_{g}$ on the first homology group of the surface $S$ and
the kernel of the representation $\rho$ is, by definition, the Torelli group
$Tor^n_{g}$. The Dehn twists $D_a$, with $a$ such that $S-a$ is not connected
("bounding curve") and the products $D_aD_b^{-1}$, with $a$ and $b$ not disconnecting
$S$ but such that $S-\{a,b\}$ is not connected ("bounding pair") generate the Torelli
group $Tor_g$ for $g\geq 3$.

Moreover fix an embedding of $S$ in $\mathbf{E}^3$ as in the
figure below:
\begin{figure}
[tbh]
\begin{center}
\includegraphics[
height=3.5147cm, width=8.843cm
]%
{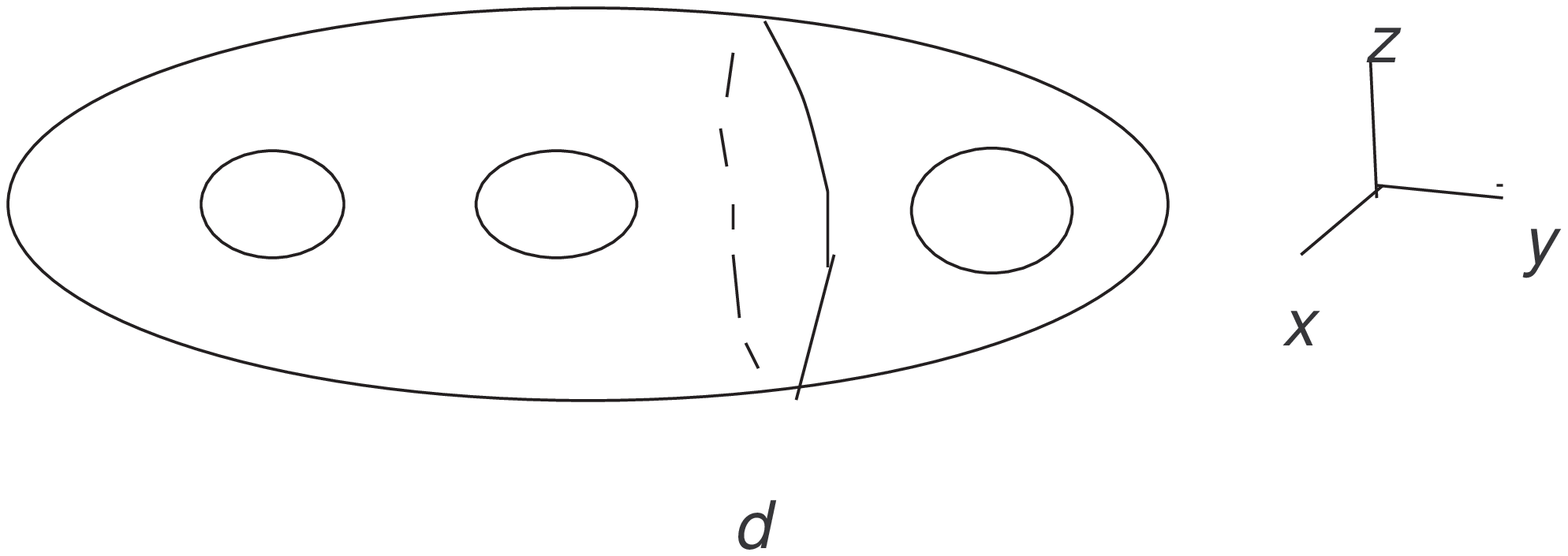}%
\caption{}%
\label{Fig.1}%
\end{center}
\end{figure}

The rotation of 180 degrees around the $y$-axis defines a
topological involution $i$ (called hyperelliptic involution)  on
$S$ with quotient homeomorphic to the 2-sphere $S^2$. Let
$$h_S:S\rightarrow S^2$$ be the corresponding 2 to 1 map with
$2g+2$ fixed points.

Denote by $\G_g^{H}$ the subgroup of $\G_g$ generated by those
elements which can be represented by fiber preserving
diffeomorphisms with respect to $h_S$. The group $\G_g^{H}$ is
isomorphic to the group of fiber preserving diffeomorphisms of $S$
modulo fiber preserving isotopies (for example by Th.1 of
\cite{BH}). Moreover a diffeomorphism of $S$ that is  isotopic to
identity and preserves the fibers of $h_S$ is  isotopic to
identity through  fiber preserving diffeomorphisms (cf. Th.2 of
\cite{BH}). From this it follows that $\G^H_g$ is an extension of
the mapping class group of $S^2$ with $2g+2$ marked points by
${\bf Z}/2{\bf Z}$. Set: $$Tor^H_g:=Tor_g\cap \G^H_g$$

\subsection{Moduli of curves}\label{teich} The mapping
class group is related to the moduli space of curves
 via Teichmuller theory. The Teichmuller space is the space of complex structures on
 $S$ up to isotopies that fix $\{x_1,...,x_n\}$ pointwise.
 It is a contractible complex manifold of dimension
$3g-3+n$ on which $\G^n_g$ acts properly discontinuosly with
quotient
 analytically
isomorphic to the moduli space $\mathcal{M}^n_g$ of $n$-marked
smooth projective curves of genus $g$. The action of $Tor^n_g$ is
free and the quotient $\mathcal{T}^n_g:=Tor^n_g\backslash
\mathcal{X}^n_g$, called the Torelli space, is
the moduli space of $n$-pointed smooth projective curves $C$ of
genus $g$ with a fixed symplectic basis of homology.

There are the natural projections: $$q:\mathcal{X}^n_g
\longrightarrow \mathcal{T}^n_g,\quad
p_T:\mathcal{T}^n_g\longrightarrow  \mathcal{M}^n_g, \quad
p_X=p_T\circ q:\mathcal{X}^n_g \longrightarrow \mathcal{M}^n_g.$$

\subsection{Moduli of hyperelliptic curves}
Let $H_g\subset \mathcal{M}_g$ be the moduli space of
hyperelliptic curves of genus $g$. Let $\tH$ be a connected
component of $p^{-1}(H_g)$ and $H:=q(\tH)$ the corresponding
connected component of $p_T^{-1}(H_g)$. They are complex
submanifolds of dimension $2g-1$ of $\mathcal{X}_g$ and
$\mathcal{T}_g$ respectively. The group $\G_g^{H}$ is the orbifold
fundamental group of the hyperelliptic locus $H_g$ and $Tor^H_g$
is the fundamental group of $H$.
 Let $$\pi:{\mathcal
C}\longrightarrow H,\quad \mathrm{and}\quad \pi_{\mathcal{J}}:
\mathcal{J}\mathcal{C}\rightarrow
 H$$ be the universal families of curves and jacobians on $H$.

\begin{lem} The family  $\pi:{\mathcal
C}\longrightarrow H$
 has $2g+2$ sections $\tilde{q}_i:H\longrightarrow \mathcal{C}$
corresponding to the Weierstra{\ss} points sets. In particular we
ask the first two
 sections to satisfy
 $h_t(\tilde{q}_1(t))=0$ and $h_t(\tilde{q}_2(t))=\infty$ and we denote the third one
  by $\tilde{p}$.
\label{wp}
\end{lem}

\begin{proof}
The set $W_t$ of Weierstra{\ss} points of $C_t:=\pi^{-1}(t)$ defines a covering
$\mathcal{W}\rightarrow H$ of degree $2g+2$ . The lemma states that $\mathcal{W}$ has
$2g+2$ connected components defining the sections $\tilde{q}_i$. Consider the
universal family of theta characteristics: $$\mathcal{S}_g:=\{L\in
Pic(\mathcal{C}_{\mathcal{T}}/T_g),\;L^2\simeq
\omega_{\pi_{\mathcal{T}}}\}\longrightarrow \mathcal{T}_g.$$ The fiber $S_t$ over
$t\in \mathcal{T}_g$ is the space of the theta-characteristics over $C_t$. The set
$S_t$ can be identified with the set $Q$ of all the quadratic forms
 $H_1(C_t,{\bf Z}/2{\bf Z})\rightarrow {\bf Z}/2{\bf Z}$ compatible with the
intersection product. Since $Tor_g$ acts trivially on the first
homology group, we have that:
$$\lambda:\mathcal{S}_g\longrightarrow  Q\times \mathcal{T}_g$$ is
an isomorphism. From the natural inclusion $${\mathcal
W}\hookrightarrow \mathcal{S}_{g|H},\quad \tilde{q}_i(t)\mapsto
\mathcal{O}((g-1)\tilde{q}_i(t))$$ and  the isomorphism $\lambda$,
it follows that $\mathcal{W}$ is the union of $2g+2$ connected
components which define the sections.

\end{proof}
\subsection{The normal function extending $reg(Z)$}
 Consider on $H$ the
variations of Hodge structures $R^2\pi_{\mathcal{J}*}\mathbf{Z}\simeq
\wedge^2\mathcal{R}^{1}\pi_{*}\mathbf{Z}$ and $\mathcal{P}_2\simeq
ker(\wedge^2\mathcal{R}^{1}\pi_{*}\mathbf{Z}\rightarrow
\mathcal{R}^{2}\pi_{*}\mathbf{Z})$ which extend the second and second primitive
cohomology of the jacobian. Associated to them there are the families of intermediate
jacobians:
 $$\mathcal{J}_2 :=
(\mathcal{F}^1\mathcal{R}^{2}\pi_{\mathcal{J}*}\mathbf{C})^*/
\mathcal{R}^{2}\pi_{\mathcal{J}*}\mathbf{Z}^*\rightarrow H, \quad \quad
\mathcal{J}_{2prim}:= (\mathcal{F}^1\mathcal{P}_{2\mathbf{C}})^*/
\mathcal{P}_2^*\rightarrow H,$$
 extending $J_2(JC)$ and $J_2(JC)_{prim}$ (cf. \ref{def1reg}, \ref{def2reg}).

Since $Tor_g$ acts trivially on $H^1$ the local system
$\mathcal{R}^{1}\pi_{*}\mathbf{Z}$, the fibrations $\mathcal{J}\mathcal{C}$,
$\mathcal{J}_2$ and $\mathcal{J}_{2prim}$ are topologically trivial on $H$.
 Let
$$p_{J_2}:\mathcal{J}_2\rightarrow J_2(JC),\quad\quad
p_{J_{2prim}}:\mathcal{J}_{2prim}\rightarrow J_2(JC)_{prim}$$ be
the projections onto the fibers.

 By Lemma \ref{wp} $H$ is a fine moduli space for hyperelliptic
curves with marked Weierestra{\ss} points, hence the Collino cycle
extends to a family of higher cycles $\mathcal{Z}$ on the
associated family of jacobians.
 The construction of the regulator images of any fiber $Z_t$ of
 $\mathcal{Z}$
extends to  normal functions $\mathrm{R}_{\mathcal{Z}}$ and
$r_{\mathcal{Z}}$ i.e. to holomorphic sections of $\mathcal{J}_2 $
and $\mathcal{J}_{2prim}$: $$\begin{array}{ccc} &&\mathcal{J}_2\\
&\stackrel{\mathrm{R}_{\mathcal{Z}}} \nearrow&\downarrow\\
H&\stackrel{id}\longrightarrow &H
\end{array}\quad \begin{array}{ccc}
&&\mathcal{J}_{2prim}\\ &\stackrel{r_{\mathcal{Z}}}
\nearrow&\downarrow\\ H&\stackrel{id}\longrightarrow &H
\end{array}.$$

\subsection{Remark} The construction of these normal functions
could have been done (as Collino does) on a finite covering of
$H_g$ given for example by the moduli space of hyperelliptic
curves with a convenient level structure. The reason why we  work
on $H$ is that, using the projections to the fibers $p_{J_{2}}$
and $p_{J_{2prim}}$, we can forget about the $Sp_g$-contribution
to monodromy as in the following.

 \subsection{The induced homomorphism}  We are interested
in  the homomorphism of fundamental groups induced by the compositions $p_{J_2}\circ
\mathrm{R}_{\mathcal{Z}}$ and $p_{J_{2prim}}\circ r_{\mathcal{Z}}$. To prove that
$reg(Z)$ is not zero for a general hyperelliptic it is enough to prove that
$$(p_{J_{2prim}}\circ r_{\mathcal{Z}})_*:\pi_1(H)=Tor_g^H\rightarrow H_2(JC)_{prim}$$
is not trivial.

In order to prove this we compute the image of the class of loops $\lambda_d$ in $H$
based at $[C]$, that correspond to a Dehn twist $D_d$ of a bounding curve $d$ on $C$,
invariant with respect to the hyperelliptic involution (thus $D_{d}\in Tor^H_g$). The
loop $\lambda_d$ lifts to a path $\tilde{\lambda}_d:[0,1]\rightarrow \tilde{H}$ with a
parametrization such that $\tilde{\lambda}_d(0)=[C]$ and
$\tilde{\lambda}_d(1)=[D_{d}C]$ in $\tilde{H}$. Restrict the universal family of
curves to $\tilde{\lambda}_d$: $$\begin{array}{ccccc}
\mathcal{C}_{\mathcal{X}|\tilde{\lambda}_d}&=:&\mathcal{C}_{\tilde{\lambda}_d}&\supseteq&C_t\\
&\pi_{X|\tilde{\lambda}}&\downarrow&&\downarrow\\
&&\tilde{\lambda}_d&\ni&\tilde{\lambda}_d(t).
\end{array}$$
For any $t$ consider the holomorphic map $h_t:C_t\rightarrow \mathbf{P}^1$ such that
$\tilde{q}_1(\tilde{\lambda}(t))=h^{-1}_t(0)$ and
$\tilde{q}_2(\tilde{\lambda}(t))=h^{-1}_t(\infty)$, corresponding to the topological
quotient $h_S:S\rightarrow S^2$. In particular $h_0=h$ and $h_1=h\circ D_{d}=:h_{d}$.
Set $\gamma_t:=h^{-1}_t([0,\infty])$ and $\gamma_{d}:=\gamma_1.$

The expression of $\reg(Z)$ obtained in Th.\ref{regreg} (see Remark \ref{closed})
allows us to lift the normal function $\mathrm{R}_{\mathcal{Z}}$ along
$\tilde{\lambda}_d$ to a section of
$(\mathcal{F}^1\wedge^2R^1\pi_{X|\tilde{\lambda}_d*}\mathbf{C})^*$ in the following
way. For all $t\in [0,1]$, let $\phi_t$ and $\psi_t$ be closed 1-forms on $C_t$, with
$\psi_t$ of type (1,0). It is enough to define the lifting of
$\mathrm{R}_{\mathcal{Z}|\lambda_d}$ along $\tilde{\lambda}_d$ as
 $$\tilde{\mathrm{R}}_{\tilde{\lambda}_d}:[0,1]\rightarrow
(\mathcal{F}^1\wedge^2R^1\pi_{X|\tilde{\lambda}_d*}\mathbf{C})^*$$ on the classes of
$\phi_t\wedge \psi_t$: $$\tilde{\mathrm{R}}_{\tilde{\lambda}_d}(t)([\phi_t\wedge
\psi_t])= 2\int_{C_t-\g_t}log(h_t)\phi_t\wedge \psi_t+ 2\pi i \int_{\g_t} (\phi_t
\psi_t- \psi_t \phi_t) .$$

By covering theory we have: $$(p_{J_2}\circ
\mathrm{R}_{\mathcal{Z}})_*([\lambda_d])=(1/2\pi)(
\tilde{\mathrm{R}}_{\tilde{\lambda}_d}(1)-\tilde{\mathrm{R}}_{\tilde{\lambda}_d}(0))\in
H_2(JC,\mathbf{Z}). $$

\begin{prop} Let $\lambda_d$ be a loop in $H$ with basis point $[C]$, whose homotopy class
 corresponds to the
Dehn twist of a bounding curve $d$ invariant for the hyperelliptic
involution and splitting $C$ in a component $S_1$ containing $q_1$
and a component $S_2$ containing $q_2$. Let $\phi$ and $\psi$ be
closed 1-forms on $C$, with $\psi$ of type (1,0), then
 $$(p_{J_2}\circ
\mathrm{R}_{\mathcal{Z}})_*([\lambda_d])([\phi\wedge
\psi])=4<[\phi]_{S_2},[\psi]_{S_2}>_{S_2}$$ where $<,>_{S_2}$ is the intersection form
on $S_2$.\label{loo}
\end{prop}

\begin{proof}
Notice that, since $D_{d}$ acts trivially in homology and by
Remark \ref{closed}, we can choose $\phi_1=\phi_0:=\phi$ and
$\psi_1=\psi_0:=\psi$. Thus
$$\tilde{\mathrm{R}}_{\tilde{\lambda}_d}(1)-\tilde{\mathrm{R}}_{\tilde{\lambda}_d}(0)(\phi\wedge
\psi)= 2\left(\int_{C_d-\g_d}\log(h_d)\phi\wedge
\psi-\int_{C-\g}\log(h)\phi\wedge \psi\right)$$ $$+ 2\pi i\left(
\int_{\g_d} (\phi \psi- \psi \phi)-\int_{\g} (\phi \psi- \psi
\phi)\right).$$ The proposition follows then from the following
two equalities:

\begin{equation}\label{monreg1}\int
_{C_d-\g_d}\log(h_d)\phi\wedge\psi-\int_{C-\g}\log(h)\phi\wedge\psi
=4\pi i <[\phi]_{S_2},[\psi]_{S_2}>_{S_2}-4\pi i
\int_{d}\phi\psi.\end{equation}
\begin{equation}\label{monreg2}
\int_{\gamma_d}(\phi\psi-\psi\phi)-\int_{\gamma}(\phi\psi-\psi\phi)=4
\int_d\phi\psi.\end{equation} First we prove \ref{monreg1}. By the assumptions, the
separating curve $d$ is the preimage of a curve $d_0$ on $\mathbf{P}^1$ whose
complement is the union of two open disks that are neighborhoods of 0 and $\infty$
respectively. The surface with boundary $S_1$ is the preimage of a disk containing 0
and $S_2$ is the preimage of a disk containing $\infty$. We assume that $d_0$
intersects $[0,\infty]$ transversely in a single point. Thus the Dehn twist $D_{d_0}$
carries $[0,\infty]$ to $[0,\infty] +d_0$. The square $D^2_{d_0}$ lifts to the Dehn
twist $D_{d}$. Its effect on the integral is via the multivaluedness of the logarithm:
this will change  by $4\pi i$ on $S_2$. Then the difference between the 2 integrals is
$4\pi i\int_{S_2}\phi\wedge \psi$. To calculate this, first observe that since the
boundary curve $d$ is null homologous, the restriction of $\phi$ to $d$ is exact. Let
$\rho$ be a smooth function $\rho$ on $C$ such that $\phi_0=\phi-d\rho$ vanishes on a
neighborhood of $d$ and such that $\rho$ is zero in the point $d\cap \gamma^+$ .
Clearly $\phi_0$ defines the same homology class as $\phi$. Its restriction to $S_2$
vanishes near the boundary so, by the De Rham theorem: $$\int_{S_2}\phi_0\wedge\psi
=<[\phi]_{S_2},[\psi]_{S_2}>_{S_2}.$$ It remains to compute
$$\int_{S_2}d\rho\wedge\psi.$$ If $d$ is orientated as the boundary of $S_1$, then
Stokes' theorem implies that this is equal to $-\int_{d}\rho\psi$, which is, by
definition the iterated integral and Lem.\ref{ii}.(3), $-\int_{d}\phi\psi$.

We prove now equality \ref{monreg2}. The Dehn twist $D_{d}$
carries $\gamma$ to $\gamma_{d}:=\gamma+2 d$, with the chosen
orientation of $d$. Then, by Lem.\ref{ii}.(1) and (2) and the fact
that $\int_{d}\phi=\int_{d}\psi=0$: $$
\int_{\gamma_{d}}(\phi\psi-\psi\phi)-\int_{\gamma}(\phi\psi-\psi\phi)=2
\int_{d}(\phi\psi-\psi\phi)=4\int_{d}\phi\psi.$$

\end{proof}

\begin{cor} Keep the notation of Prop.\ref{loo}. Choose a symplectic
basis $\{A_k,A_{g+k}\}_{k\leq g}$ of $H_1(C,{\bf Z})$  such that $\{A_k,
A_{g+k}\}_{k\leq g_1}$ is a symplec\-tic basis for $H_1(S_1,{\bf Z})$ and
$\{A_k,A_{g+k}\}_{g_1<k\leq g}$ is a symplectic basis for $H_1(S_2,{\bf Z})$, then
$$(p_{J_2}\circ \mathrm{R}_{\mathcal{Z}})_*([\lambda_d])= 4\sum_{k=g_1+1}^g A_k\wedge
A_{g+k}.$$\label{corloo}
\end{cor}

\begin{proof}
The intersection form on $S_2$, seen as an element in $\wedge^2H^1(C,\mathbf{Z})$ is,
up to constant: $$<,>_{S_2}=\sum_{k=g_1+1}^g A_k\wedge A_{g+k}$$
\end{proof}

\subsection{Remark} Notice that
$H^2(JC,\mathbf{Z})_{prim}^*\simeq H_2(JC,\mathbf{Z})_{prim}\simeq
\wedge^2H_1(C,\mathbf{Z})/<\omega>$, where $\omega:=\sum_{k=1}^g A_k\wedge A_{g+k}$ is
the dual of the polarization $\Omega$.
\begin{cor}
With the notation of Cor.\ref{corloo}, we have: $$(p_{J_{2prim}}\circ
r_{\mathcal{Z}})_*([\lambda_d])= 4\sum_{k=g_1+1}^g A_k\wedge
A_{g+k}\quad\mathrm{mod}\quad <\omega>.$$In particular $ (p_{J_{2prim}}\circ
r_{\mathcal{Z}})_*$ is not trivial and hence $reg(Z)$ is not zero for general
hyperelliptic.\label{torsione}
\end{cor}
\begin{proof}
It follows directly from the previous corollary, since
$r_{\mathcal{Z}}$ is the composition of $\mathrm{R}_{\mathcal{Z}}$
with the natural projection $\mathcal{J}_2\rightarrow
\mathcal{J}_{2prim}$.
\end{proof}

\section{Monodromy of the extensions}

\subsection{}The goal of this section is to extend the extension classes of Sect.3,
 to normal functions on the moduli space
$H$ and to study the induced homomorphism on fundamental groups. We compare the
homomorphism induced by the normal function extending $Pe$ with the one induced by
$R_{\mathcal{Z}}$ (see Cor.\ref{uu}).
\subsection{The normal functions extending
$e_s$ and $Pe$} Any finite dimensional $Sp_g(\mathbf{Z})$ representation $V$ has a
natural Hodge structure which can be extended to a variation of Hodge structures
$\mathcal{V}$ on any fine moduli space of curves.  Moreover if $V$ has negative
weight, its intermediate jacobian $J(V)$ (see \ref{JV}) can be extended to a
corresponding intermediate jacobians fibration
 $$\mathcal{J}(\mathcal{V}):={{\mathcal{V}_\mathbf{C}}\over{\mathcal{F}^0+
 \mathcal{V}}}.$$
 We will restrict
ourselves to the case in which
 this moduli
 space is $H$.
 Since $Tor^H_g$ acts trivially on homology
 the associated bundles of intermediate jacobians are topologically trivial.
 We shall denote by
$$p_{V}:\mathcal{J}(\mathcal{V})\rightarrow J(V)$$ the projection onto the fiber.

 We consider the case $V=Hom(\otimes^3H^1,H^1)$ and denote by  $\mathcal{E}_{s}$ ($s=1,2$) and
$\mathcal{P}\mathcal{E}$ the normal functions which associate to
 the curve $C$
the extensions $e_{s}$ and $Pe$ fitting in the following commutative diagram:
\begin{equation}\label{diagram}\begin{array}{ccccccc} &&\mathcal{J}(\mathcal{V})&
\stackrel{\phi}\rightarrow
&\mathcal{J}_{2prim} &&
\\
&\stackrel{\mathcal{E}_{2}-\mathcal{E}_{1}}
\nearrow&\downarrow&&\downarrow&\stackrel{\mathcal{P}\mathcal{E}}\nwarrow&\\
H&\stackrel{id}\rightarrow&H&\stackrel{id}\rightarrow&H&\stackrel{id}\leftarrow&H.
\end{array}
\end{equation}

\subsection{Monodromy of $\mathcal{E}_s$}

The pairs of sections $(\tilde{q}_s,\tilde{p})$ ($s=1,2$), of Lemma \ref{wp}, define
two different  inclusions: $$j_{q_s,p}:H\hookrightarrow T_g^2.$$ We will show that the
homomorphism $$(p_{V}\circ\mathcal{E}_{s})_*:Tor_g^H\longrightarrow
Hom(\otimes^3H^1,H^1))\simeq Hom(H_1,\otimes^3H_1)$$ factorizes via $(j_{q_s,p})_*$
and a natural homomorphism associated to the action of the mapping class group on
$J_{q_s,p}/J^4_{q_s,p}$ described in the following.
 Let
$\pi:=\pi_1(C-\{q\},p)$ be the fundamental group of the punctured
curve $C-\{q\}$ in $p$. The mapping class group $\G^2_g$ acts
naturally on $\pi$ and on its lower central series of $\pi$. The
action on $\pi$ induces a natural action on the $J$-adic quotients
of the group algebra of the fundamental group ($J=J_{q,p}$ in the
notation of Sec.3) that we will briefly illustrate (see
\cite{J1},\cite{J2}).

Let $\{\pi^{(k)}\}$ $(k\geq 1)$ be the lower central series of
$\pi$, namely $\pi^{(1)}:=\pi$ and $\pi^{(k+1)}=[\pi^{(k)},\pi]$,
with nilpotent quotients $$N_k:=\pi/\pi^{(k)}\quad {\rm and}\quad
L_k:=\pi^{(k)}/\pi^{(k+1)}$$ fitting in the sequence
$0\longrightarrow L_k\longrightarrow N_{k+1}\longrightarrow N_k
\longrightarrow 0.$

The action of $\G^2_g$ on $\pi$ induces an action on the quotients
 $N_k$:
\begin{equation}\label{rok}
\rho_k:\G^2_g\longrightarrow Aut(N_k) .
\end{equation}
Notice that $N_2=H_1$, so $\rho_2=\rho$ and $ker\rho_2=Tor^2_g$.
The action of $\G^2_g$ on $L_k$ factorizes through $\rho$, thus
$L_k$ is a $Sp_{g}$-module. More precisely $\oplus_k L_k$ is the
free Lie algebra (over {\bf Z}).
 In particular $L_k$ is
torsion free and there are the following identifications: $$L_2\simeq
\wedge^2H_1,\quad L_3\simeq (\wedge^2H_1\otimes H_1)/\wedge^3H_1.$$ Johnson's
homomorphisms $\tau_k$ in the punctured case are:

\begin{equation}\label{tauk}
\tau_k:ker\rho_k\longrightarrow Hom(H_1,L_k)\quad f\mapsto
[\overline {n}\mapsto f(n)n^{-1}\;\mathrm{mod}\;\pi^{(k+1)}]
\end{equation}
where $n$ is any lifting of $\overline {n}$ to $N_{k+1}$. By construction
$ker\rho_{k+1}=ker\tau_k$.

\begin{lem}\label{j}
The inclusion $j:\pi\hookrightarrow J$, $\alpha\mapsto \alpha-1$ satisfies
$j(\pi^{(k)})\subset J^k$ and induces an injective homomorphism $j_k:L_k\rightarrow
J^{k}/J^{k+1}$.\end{lem}

\begin{proof} The inclusion $j(\pi^{(k)})\subset J^k$ can be checked by
induction and it implies that $\pi^{(k)}\subset \pi^k_J$ where
$\pi^k_J:=\{\alpha\in \pi|\;\alpha-1\in J^k\}$. Moreover the group
$\pi/\pi^k_J$ has no torsion and $\pi^k_J/\pi^{(k)}$ is a torsion
group, hence
 $j_k:L_k\otimes_{{\bf Z}}{\bf Z}\hookrightarrow J^{k}/J^{k+1}
\otimes_{{\bf Z}} {\bf Z}$ and, since $L_k$ has no torsion, $j_k:L_k\hookrightarrow
J^{k}/J^{k+1}$.
\end{proof}

We define analogues of Johnson's homomorphisms associated to the action on $J$-adic
quotients:
\begin{lem} The homomorphim $$ \tau_{kJ}:ker\rho_k\longrightarrow
Hom(H_1,J^k/J^{k+1}),\quad\quad\tau_{kJ}:f\mapsto \{[\alpha]
\mapsto f(\alpha)-\alpha\;\mathrm{mod}\;J^{k+1}_p \}$$ is well
defined
 and $$\tau_{kJ}=H(j_k)\circ \tau_k,$$

 where $H(j_k):Hom(H_1,L_k)\rightarrow Hom(H_1,J^k/J^{k+1})$,
 $H(j_k):f\mapsto j_k\circ f$.
\end{lem}

\begin{proof} By Lem.\ref{j}
 $ker\rho_k\subset ker(\G^2_g
\longrightarrow Aut(\pi/\pi^k_J))=\ker(\G^2_g \longrightarrow
GL(J/J^k))$ hence $\tau_{kJ}$ is  well defined. Moreover
$\tau_{kJ}(f)(a)=j_k\tau_k(f)(a)$ where $f\in \G^1_g(k)$, $a\in
H$, because: $f(n)-n=(f(n)n^{-1}-1)n=f(n)n^{-1}-1\;{\rm
mod}\;J^{k+1}$.

\end{proof}

\begin{lem} The homomorphism $(j_{q_s,p})_*$ satisfies:
\begin{equation}
(j_{q_s,p})_*(Tor^H_g)\subset ker\rho_3=ker\tau_2.\end{equation} \label{Mor}
\end{lem}

\begin{proof}
 Since the hyperelliptic involution acts as $-I$ on $H_1$, it also acts as $-I$ on
$Hom(H_1,L_2)\simeq Hom(H_1,\wedge^2H_1)$  while it acts
trivially on $(j_{q_s,p})_*(Tor^H_g)$.

\end{proof}

\begin{lem} The following equality holds:
\begin{equation}
(p_{V}\circ\mathcal{E}_{s})_*=\tau_{3J}\circ (j_{q_s,p})_*.\end{equation}
\label{sofia}
\end{lem}

\begin{proof}
 The homomorphism
$(p_{V}\circ\mathcal{E}_{s})_*$ is determined by the action of
$(j_{q_s,p})_*(Tor_g^H)\subset \G^2_g$ on $J_{q_s,p}/J^4_{q_s,p}$.
Since, by the previous lemma, this action is trivial on $N_3$, it
defines in fact $\tau_{3J}$.

\end{proof}

\subsection{Monodromy of $\mathcal{P}\mathcal{E}$} We want to compute
 $(p_{J_2prim}\circ\mathcal{P}\mathcal{E})_*$ on the Dehn twist $D_d$ in $Tor^H_g$ of a simple curve
$d$ separating $q_1$ and $q_2$ exactly as in Sect.4. We first compute $(p_V\circ
(\mathcal{E}_2-\mathcal{E}_1))_*$:

\begin{prop}
\label{altrononzero} Let $D_d\in Tor^H_g$ be the Dehn twist of the simple curve $d$
separating $q_1$ and $q_2$ and such that $p$ is in the same component of $q_1$. It
holds
$$(p_{V}\circ(\mathcal{E}_{1}-\mathcal{E}_{2}))_*(D_d)=\left\{\begin{array}{ccccc} \;
\omega\wedge A_k &{\rm if}& g_1+1\leq l\leq g &{\rm or}& g+g_1+1 \leq l\leq 2g\\
\;0&&{\rm otherwise}&&
\end{array}\right.
$$ where $\omega:=\sum_{k=1}^g A_k\wedge A_{g+k}$ is the dual of the polarization
$\Omega$.
\end{prop}

\begin{proof}
 By Lem.\ref{sofia}
 $$(p_{V}\circ(\mathcal{E}_{1}-\mathcal{E}_{2}))_*(D_d)=\tau_{3J}
 (j_{q_2,p})_*(D_d)(j_{q_1,p})_*(D_d)^{-1}).$$
Fixing the isomorphism $\G^2_g\simeq Aut^+(\pi_1(C-\{q\},p))$, the element
$(j_{q_2,p})_*(D_d)(j_{q_1,p})_*(D_d)^{-1})\in Tor^2_g$ can be identified with the
element $D_{d'}D_{d"}^{-1}$ as in Fig.2, where $d'$ and $d"$ are homotopic to $d$ and
bound a cilinder containing the missing point $q$.
\begin{figure}
[tbh]
\begin{center}
\includegraphics[
height=3.5147cm, width=8.843cm
]%
{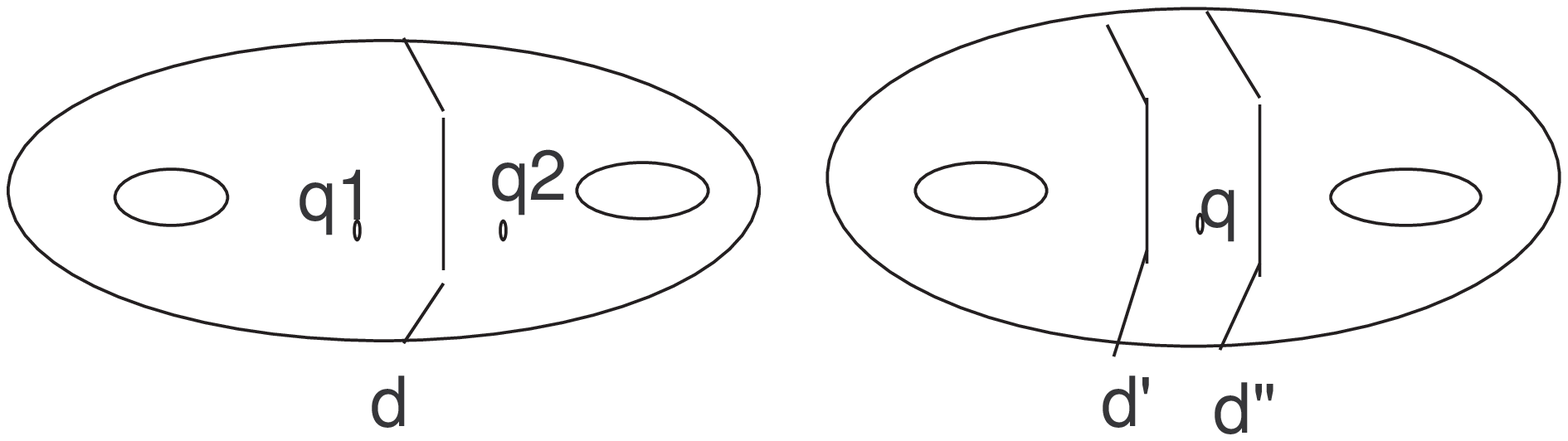}%
\caption{}%
\label{Fig.2}%
\end{center}
\end{figure}

The loops $\al_l$, for $1\leq l\leq g_1$ and $g+1\leq l\leq g+g_1$, can be chosen not
to intersect $d'$ and $d"$ thus the action of $D_{d'}D_{d"}^{-1}$ on them  is the
identity, while on $\al_l$ for $g_1+1\leq l\leq g$ and $g+g_1+1 \leq l\leq 2g$ the
action is given by the conjugation by $\delta_q$ where $\delta_q$ bounds a punctured
disk around $q$ and satisfies $\delta_q=\prod_k[\al_k\al_{g+k}]$. Hence, for $1\leq
l\leq g_1$ and $g+1\leq l\leq g+g_1$ $$\tau_3(D_{d'}D_{d"}^{-1})(A_l)=0,$$ while for
$g_1+1\leq l\leq g$ and $g+g_1+1 \leq l\leq 2g$
$$\tau_3(D_{d'}D_{d"}^{-1})(A_l)=[\prod_{k=1}^g[\al_k\al_{g+k}],\al_l]\in L_3.$$ Via
the identification $L_3\simeq(\wedge^2H_1\otimes H_1)/\wedge^3H_1$ and the inclusion
$J_3:L_3\hookrightarrow J^3/J^4=\otimes^3H_1$ we get
$$\tau_{3J}(D_{d'}D_{d"}^{-1})(A_l)=\left(\sum_{k=1}^g A_k\wedge A_{g+k}\right)\wedge
A_l=\omega\wedge A_l.$$
\end{proof}
\begin{cor}\label{uu}
$$(p_{J_{2prim}}\circ\mathcal{P}\mathcal{E})_*(D_d)=2(2g+1)\sum_{k=g_1+1}^{g}A_k\wedge
A_{g+k}\; mod\; \omega =(2g+1)(p_{J_{2prim}}\circ r_{\mathcal{Z}})_*(D_d).$$
\end{cor}
\begin{proof}
The fiber bundle map $\phi:\mathcal{J}(\mathcal{V})\rightarrow \mathcal{J}_{2prim}$ of
diagram \ref{diagram},  restricted to a fiber, is the surjective
 homomorphism $\Phi:J(V)\simeq Ext_{MHS}(\otimes^3H^1,H^1)\rightarrow J_2(JC)_{prim}\simeq
 Ext_{MHS}(\wedge^2H^1,\mathbf{Z})$ described in Sect.\ref{ext}: i.e. it is obtained by pulling back along
 $J_{\Omega}$, tensoring by $H^1$ on the left, pushing down along $\Pi$ and pulling
 back along the inclusion $\wedge^2H^1_{prim}\hookrightarrow \otimes^2H^1$. The
 induced homomorphism between the fundamental groups of the two intermediate jacobians
  is then: $$\Phi_*:Hom(H_1,\otimes^3H_1)\rightarrow \wedge^2H_{1prim}$$ $$f\mapsto
\sum_kA_k\wedge J^t_{\Omega}f(A_{g+k})-A_{g+k}\wedge
J^t_{\Omega}f(A_{k})\:mod\;\omega$$ where $J^t_{\Omega}:\otimes^3H_1\rightarrow
H_1,\;A_l\otimes A_m\otimes A_n\mapsto
\sum_k(\delta_{lk}\delta_{m(g+k)}-\delta_{l(g+k)}\delta_{mk})A_n$ is the dual of
tensoring by $\Omega$ from the left (cf. \ref{Omega}). The first equality then follows
by direct computation using the equality: $$p_{J_{2prim}*}\circ \phi_*=\Phi_*\circ
p_{V*},$$ Prop.\ref{altrononzero} and the explicit expression of $\Phi_*$. The second
equality is the statement of Cor.\ref{corloo}.

\end{proof}

\section{Final Remarks}
\subsection{An alternative proof of Th.\ref{risultato} } The computations of Cor.\ref{uu}  and Cor.\ref{corloo} can be
easily improved to show that two normal functions $(2g+1)r_\mathcal{{Z}}$ and
$\mathcal{P}\mathcal{E}$ induce the same homomorphism on fundamental groups, without
applying Th.\ref{risultato}. In fact this last monodromy computation can be used to
give
 an alternative proof of  Th.\ref{risultato}, following the steps in which  Hain
 proved in a different way in \cite{HD4} the result of Hain-Pulte
on $C_p-C_p^-$ and the extension $J_p/J_p^3$. To apply the rigidity argument of
\cite{HD4}Cor.6.4 to our case one just needs  to show the following two facts. First
of all both $(2g+1)r_\mathcal{{Z}}$ and $\mathcal{P}\mathcal{E}$ have to define a good
VMHS in the sense of Saito \cite{Sa}. Secondly one has to exhibit a hyperelliptic
curve $C$ with  Weierstra{\ss} points $p$, $q_1$ and $q_2$, for which
$(2g+1)reg(Z)=Pe$. This last point could be achieved by taking a hyperelliptic curve
with a further involution exchanging $q_1$ and $q_2$ and fixing $p$. For such a curve
$reg(Z)=0$ while $Pe=-Pe$ so $Pe$ is 2-torsion. Then one needs to show that $Pe$ is
zero, for example by showing that there are no sections of
$\mathcal{J}_{2prim}\rightarrow H$ of order 2. The variation defined by
$\mathcal{P}\mathcal{E}$ is good since the set of the $J/J^k$ form a good VMHS over
any fine pointed moduli space of curves (\cite{HD2}). About $r_\mathcal{{Z}}$ we can
argue as follows. Let $W=\sum_i(V_i,f_i)$ be a cycle in $CH^n(X,1)$, with $X$ smooth
projective of dimension $n$. The image of $\reg(W)$ under
 the map: $J_2(X)(=Ext_{MHS}(\mathbf{Z}(n-1),H^{2n-2}(X))\rightarrow
 Ext_{MHS}(\mathbf{Z}(n-1),H^{2n-2}(X)/<[V_i]>)$ is described by the
 following construction (see for example \cite{Voisin} (0.5) or
\cite{Muller}).
 Set $U=X-|W|$ and look at the long exact sequence:
 $$...\rightarrow H^{2n-2}_{|W|}(X)\rightarrow
 H^{2n-2}(X)\rightarrow H^{2n-2}(U)\rightarrow...$$ The cycle $W$ gives rise
 to a class in $F^n\cap H^{2n-2}_{|W|}(X, \mathbf{Z})$ hence by pullback it gives rise to an
 extension

 $$0\rightarrow H^{2n-2}(X)/<[V_i]>\rightarrow E\rightarrow \mathbf{Z}(n-1)\rightarrow 0.
 $$
The class of this extension is exactly the image of $\reg(W)$.
 Since $E$ is a sub MHS of the cohomology of the open variety $U$,  when we extend
 the construction to families we get a "geometric"  VMHS,  hence it is admissible for Steenbrink and Zucker
 (cf.\cite{SZ} and \cite{Al} and
then good for Saito. In the case of the Collino cycle $Z=(C_1,h_1)+(C_2,h_2)$,
$H^{2g-2}(JC)/<[C_i]>=H^{2g-2}(JC)/<[C]>\simeq H^{2}(JC)_{prim}^*$, hence
$Ext_{MHS}(\mathbf{Z}(g-1),H^{2g-2}(JC)/<[C_i]>)=J(C)_{2prim}$ and the class of $E$ is
$reg(Z)$.

\subsection{Degeneration} As last remark we illustrate a heuristic argument
about the monodromy of the normal function associated to $reg(Z)$. This argument leads
to the formula Cor.\ref{corloo} without writing $reg(Z)$ in terms of integrals on the
curve $C$ as in Th.\ref{reg}. The starting point is the idea, explained in Collino
(\cite{Col}1.1), of viewing the higher cycle $Z$ as a degeneration of the Ceresa
cycle, when two points are identified. More precisely, set $ \tilde{C}_1:=\psi_1(C)$
and $ \tilde{C}_2:=\psi_2(C)$ where $$\psi_1:C\longrightarrow JC\times \mathbf{P}^1
,\quad x\mapsto (x-q_1,h(x)) \quad$$ $$\psi_2:C\longrightarrow JC\times
\mathbf{P}^1,\quad x\mapsto (x-q_2,1/h(x))$$ and
$$\ti{Z}:=\ti{C_1}-\ti{C_2}\sim_{hom}0.$$  Denote by $\ti{\G}$ a 3-chain such that
$\pa \ti{\G}=\ti{Z}$ and let $\al \in F^1H^2(JC)$. From Stokes' theorem it follows
that: $$\int_{\ti{\G}}p_1^*(\al)\wedge p_2^*(dz/z)\;mod\;H_2(JC,\mathbf{Z})
=\mathrm{r}eg(Z)(\al) $$ where $p_1$ and $p_2$ are the projections to $JC$ and to
$\mathbf{P}^1$.

Let $\overline{M_{g}}$ be the Deligne-Mumford compactification of the moduli space of
smooth projective genus $g$ curves. Let $M^2_g\longrightarrow \overline{M_{g+1}}$ be
the map that identifies the 2 marked points, hence it associated to a smooth curve $C$
of genus $g$ an irreducible nodal curve $C'$ of geometric genus $g+1$. Topologically
this curve can be seen as the curve obtained shrinking a loop $d_0$ of a smooth curve
$G$ of genus $g+1$ to a point, as indicated in Fig.3:
\begin{figure}
[tbh]
\begin{center}
\includegraphics[
height=3.5147cm, width=8.843cm
]%
{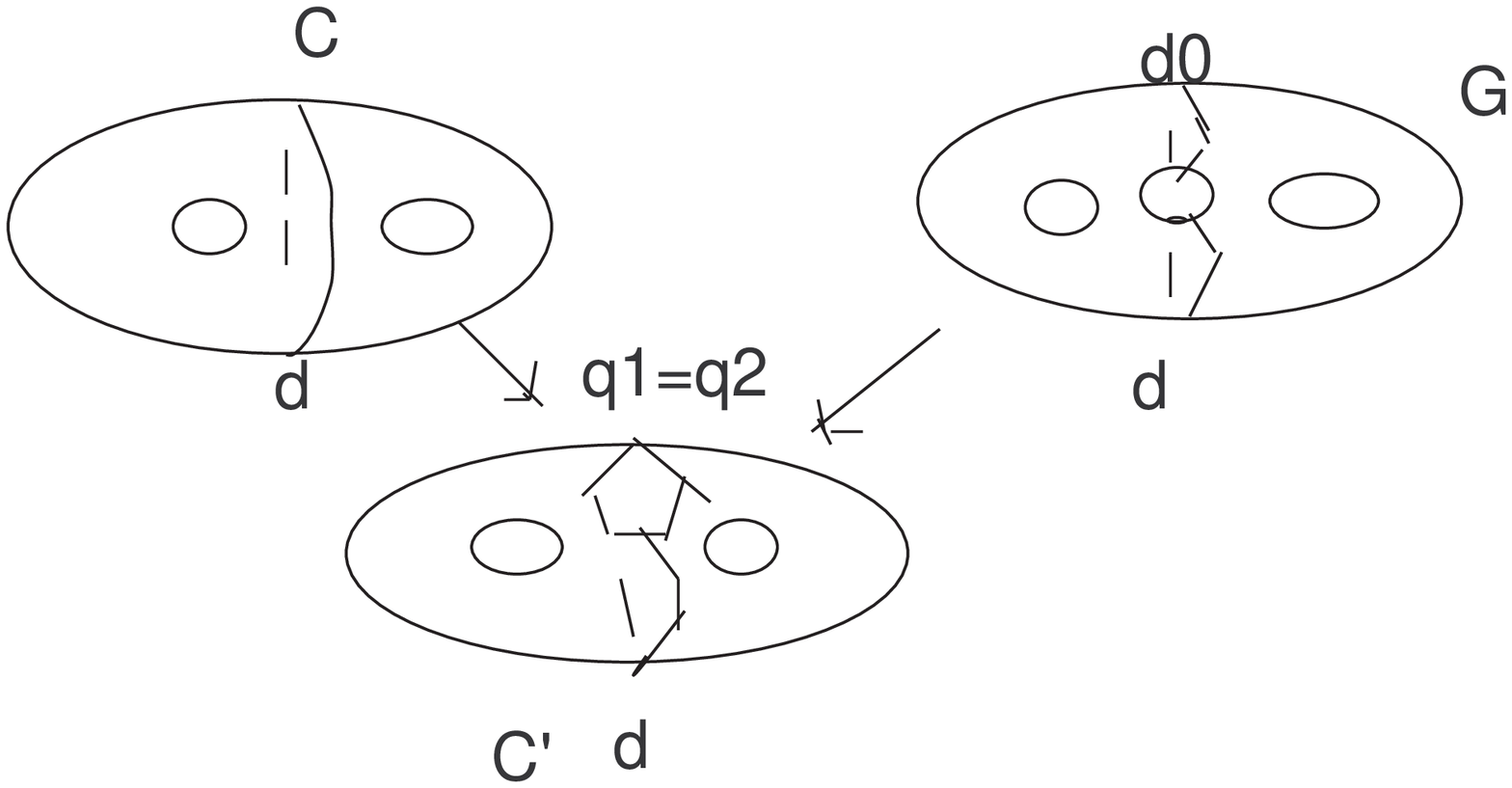}%
\caption{}%
\label{Fig.3}%
\end{center}
\end{figure}
For $C$ hyperellitic the generalized jacobian of $C'$  is (up to 2-torsion) $JC\times
\mathbf{C}^*$ so $JC\times \mathbf{P}^1$ is its compactification.
 The cycle $\ti{Z}$ can be seen as limit of the cycle $G_p-G_p^-$. Associated to it we
 have the Abel Jacobi image $Ab(G_p-G_p^-):=\int_{\G}\in
 J_3(JG):=F^2H^3(JG)^*/H_3(JG,\mathbf{Z})$, where
 $\pa \G=G_p-G_p^-$. Extending the construction
 to the universal family of curves over $T^1_{g+1}$, $Ab(G_p-G_p^-)$ extends to
 a
 normal function $\eta_{G_p-G_p^-}$ whose
 monodromy has been already determined by Hain. Let $\tau_2$ be the Johnson
 homomorphism in the
 case of a closed curve of genus $g+1$ ($\pi=\pi_1(G,p)$), Hain proved that:

$$(\eta_{G_p-G_p^-})_*(\lambda)=2\tau_2(\lambda), \quad
\forall\lambda \in Tor_{g+1}^1$$ (see Theorem 5.1 of
 \cite{HD3}). Fix a system of generators
 $\{\al_k,\al_0,\beta_0\}_{1\leq k\leq 2g}$ on $\pi_1(G,p)$
 satisfying
 the relation:
 $$\prod_{k=1}^{g_1}[\al_k,\al_{g+k}][\al_0,\beta_0]\prod_{k=g_1+1}^g
 [\al_k,\al_{g+k}]=1$$
 and such that $\al_k$ corresponds to generators of $\pi_1(C')$ and
 $\pi_1(C)$ denoted by the same letter.

 By Corollary of section 4 of \cite{J1}
 it holds $$\tau_2(D_{d}D_{d_0}^{-1})=\left(\sum_{k=g_1+1}^g[\al_k]\wedge
[\al_{g+k}]\right)\wedge [\beta_0].$$  The Dehn twist $D_{d}$ of $C$ corresponds to an
analogous Dehn twist of $C'$ and hence to the image in $C'$ of the product of Dehn
twists $D_{d}D_{d_0}^{-1}$ of $G$ in the degeneration from $G$ to $C'$. Notice that
the class $[\beta_0]$ is equal to $[d_0]$ in $C'-\{q_1\}$ and
$\int_{d_0}p_2^*dz/z=4\pi i$, thus we expected
$$(p_{J_2}\circ\mathrm{R}_{\mathcal{Z}})_*(D_{d})=\frac{1}{2\pi}8\pi
i\sum_{k=g_1+1}^gA_k\wedge A_{g+k},$$ as was proved in Cor.\ref{corloo}.

\

\end{document}